\documentclass[11pt]{amsart}
 
 \usepackage{graphicx}
 \usepackage{amsfonts}
 \usepackage{amssymb}

 \swapnumbers
 
 \advance\hoffset by -0.5 truecm

 \title[Density of elliptic geodesic]{$C^2$ densely the 2-sphere 
           \\ has an elliptic closed geodesic}
  
 \author{Gonzalo Contreras}
 \thanks{Gonzalo  Contreras was partially supported by 
  CONACYT-M\'exico grant \# 28489-E}
 \address{CIMAT \\ P.O. Box 402 \\ 36.000 Guanajuato, GTO \\ Mexico}
  \email{gonzalo@cimat.mx}
 
 \author{Fernando Oliveira}
 \address{Departamento de Matematica \\ U.F.M.G \\
 Av. Ant\^onio Carlos, 6627, Caixa Postal 702  \\
   30161-970  Belo Horizonte - MG \\ Brazil}
 \email{fernando@mat.ufmg.br}

 \newtheorem{theorem}{Theorem}[section]
 \newtheorem{Theorem}[theorem]{Theorem}
 \newtheorem{proposition}[theorem]{Proposition}
 
 \newtheorem{lemma}[theorem]{Lemma}
 \newtheorem{Lemma}[theorem]{Lemma}
 \newtheorem{definition}[theorem]{Definition}

 \newtheorem*{Kupka}{Kupka-Smale theorem for geodesic flows}

 \def\co{{\mathbb C}} 
 \def\D{{\mathbb D}}

 \def\De{\Delta}
 
 \def\cF{{\mathcal F}}
 \def\tF{{\tilde{F}}}
 \def\bG{{\mathbb{G}}}
 \def\cG{{\mathcal G}}
 \def\cH{{\mathcal H}}
 \def\bJ{{\mathbb J}}
 \def\tJ{{\tilde{J}}}
 
 \def\cO{{\mathcal O}}
 \def\P{{\mathbb P}}
 \def\bP{{\mathbb P}}
 \def\cP{{\mathcal P}}
 \def\Q{{\mathbb Q}}
 \def\hQ{{\widehat{Q}}}
 \def\cR{{\mathcal R}}
 \def\re{{\mathbb R}}
 \def\reS{{\re\times S^3}}
 \def\si{\sigma}
 \def\Si{\Sigma}
 \def\cS{{\mathcal S}}
 \def\SiGa{{\Sigma\setminus\Gamma}}
 \def\bS{{\mathbb S}}
 \def\tu{{\tilde{u}}}
 
 \def\tuas{{\tu=(a,u):\SiGa\to\reS}}
 \def\us{{u:\SiGa\to S^3}}
 \def\cU{{\mathcal U}}
 \def\cW{{\mathcal W}}
 
 \def\di{{\mathbb D}}

 \def\na{{\mathbb N}}
 \def\nulo{\emptyset}
 \def\Pg{{\Per(g)}}
 \def\oPg{{\ov{\Per(g)}}}
 \def\oPs{{\ov{\Per(\psi)}}}
 \def\cQ{{\mathcal Q}}
 
 \def\transv{\pitchfork}
 \def\ov{\overline}
 \def\scirc{{\;\scriptstyle\circ\;}}
 \def\sec{{\mathbf\Sigma}}
 \def\bS{\sec}
 \def\Y{{\mathbb Y}}
 \def\Z{{\mathbb Z}}

 \def\a{\alpha}
 
 \def\ga{\gamma}
 \def\de{\delta}
 \def\De{\Delta}
 \def\e{\varepsilon}
 \def\ep{\epsilon}
 \def\ga{\gamma}
 \def\Ga{\Gamma}
 \def\la{\lambda}
 \def\La{\Lambda}
 \def\hLa{{\widehat{\Lambda}}}
 \def\Si{\Sigma}
 \def\si{\sigma}
 \def\vr{\varphi}
 \def\om{\omega}

 \def\refteo{A }

 \DeclareMathOperator\Per{Per}
 
 \DeclareMathOperator\inte{int}
 
 \DeclareMathOperator\Sp{Sp}
 \DeclareMathOperator\length{length}
 \DeclareMathOperator\Ind{Ind}
 \DeclareMathOperator\diam{diam}
 
\begin{document}

 \begin{abstract}
 We prove that a riemannian metric on the 2-sphere or the projective plane
 can be $C^2$ approximated by a $C^\infty$ metric whose geodesic flow has an
 elliptic closed geodesic. 
 \end{abstract}
 
 \maketitle

 \bigskip

 In this paper we show how to overcome a difficulty presented by
 Henri  Poincar\'e in the $C^2$ generic case.
 In 1905, H. Poincar\'e~\cite[p. 259]{Poincare} claimed that any convex surface 
 in $\re^3$ should have an elliptic or degenerate non self-intersecting 
 closed geodesic. This is, the linearized Poincar\'e map of the geodesic
 flow at the closed geodesic has an eigenvalue of modulus 1.
 In 1980, A. Grjuntal~\cite{Grju} showed a counterexample to
 Poincar\'e's claim. Victor Donnay~\cite{Do2}, \cite{Do3} constructs an example
 of a $C^\infty$ riemannian metric on the 2-sphere $S^2$
 which has positive metric entropy and all whose closed geodesics
 but a finite number (which are degenerate) 
 are hyperbolic. Donnay's theorem is not known in positive curvature.
 It is also not known if there exists a $C^\infty$ riemannian metric
 on $S^2$ all of whose closed geodesics are hyperbolic.

 Here we prove

 \medskip

 \noindent{\bf Theorem A.}\label{teo}\quad

 {\it A riemannian metric on the 2-sphere or the projective plane
 can be $C^2$ approximated by a $C^\infty$ metric with an elliptic closed
 geodesic.}

 \medskip

 In August 2000 at Rio de Janeiro's dynamical systems 
 conference,  Michel Herman~\cite{her-an} conjectured the above result
 and  announced a proof of it in the case of positive curvature.
 His proof used, as we shall also do, Ricardo Ma\~n\'e's theory on 
 dominated splittings adapted to the geodesic flow and also
 an equivariant version of Brouwer's translation theorem. 
 We shall see that in the positive curvature case
 one can replace the use of Brouwer's translation theorem 
 by the intermediate value theorem on the interval.
 For the non convex case we use Hofer, Wysocki, Zehnder
 theory on  Reeb flows for generic tight contact forms on the
 3-sphere $S^3$.

 Theorem~\refteo is a version for geodesic flows of a theorem by
 Sheldon Newhouse~\cite{N1}. In 1977, Newhouse proved that
 if $H$ is a smooth hamiltonian on a symplectic manifold,
 0 is a regular value for $H$ and the energy level $H^{-1}\{0\}$ is
 compact, then there is a $C^2$ perturbation $H_1$ of $H$ such that
 the hamiltonian flow on $H_1^{-1}\{0\}$ is either
 Anosov or it has an elliptic closed orbit. 
 But Newhouse's arguments heavily rely on a $C^2$ closing lemma
 for hamiltonian systems, which is not known for geodesic flows.
 Newhouse's theorem applied to the hamiltonian which corresponds 
 to the geodesic flow would give a $C^2$ approximation to the metric $g$
 by a {\it Finsler} metric and not a riemannian metric.
 This is because all the known proofs for the closing lemma
 rely on local perturbations of a vector field in its
 phase space. A  perturbation of a riemannian metric is never
 a local perturbation of the geodesic flow on the unit tangent bundle
 because the support of the perturbation is a union of fibers.
 
 By a theorem of Wilhelm Klingenberg and Floris Takens~\cite{KT} 
 (see also the Kupka-Samle theorem in page~\pageref{KS})
 for any $r>0$
 one can perturb in the $C^\infty$ topology
 a degenerate or elliptic periodic geodesic 
 to make the $r$-jet of its Poincar\'e map 
 generic, in particular, with torsion.

 By Jurgen Moser's invariant circle theorem~\cite[2.11]{Moser} 
 and Birkhoff's normal form,  a $C^4$ generic elliptic
 closed geodesic in a $C^\infty$ surface has an 
 invariant torus separating the phase space. 
 Theorem~\refteo then implies that there is a 
 $C^2$ dense set of riemannian metrics on $S^2$ or $\bP^2$
 whose geodesic flow is not ergodic for the Liouville's measure. 
 The eigenvalues of the linearized
 Poincar\'e map of a generic $C^4$ elliptic geodesic is invariant
 under topological equivalences, because they can be seen as limiting
 rotation numbers on invariant torii converging to the geodesic.
 Then theorem~\refteo implies that there are no structurally stable
 geodesic flows on $S^2$ or $\bP^2$.

 Theorem~A also allows to partially generalize a result of
 Vladimir Lazutkin~\cite{Lazutkin} which says that the billiard
 map in the interior of a $C^\infty$ embedded curve in
 $\re^2$ with positive curvature is not ergodic.
 In our case we obtain that a strictly convex 
 domain in $\re^3$ with $C^\infty$ boundary in a residual
 set in the $C^2$ topology, the billiard map in its
 interior has a set of positive measure of invariant
 quasi-periodic tori. Indeed, using the KAM theorem,
 Svanidze~\cite{Svanidze} announced in $\re^3$,
 and   Valery Kovachev and Georgi Popov~\cite{KP2, KP1} 
 proved this result in $\re^n$, $n\ge 3$, provided that 
 the geodesic flow  on the boundary has an elliptic 
 periodic geodesic which is $k$-elementary, $k\ge 5$.

 \bigskip

 We sketch the proof of theorem~\refteo.
 Since $S^2$ is a double cover of $\bP^2$,
 the result on $\bP^2$ can be inferred from the
 result on $S^2$.
 Let $\cH(S^2)$ be the set of $C^\infty$ riemannian
 metrics on $S^2$ all of whose closed geodesics are hyperbolic.
 Let $\cF^1(S^2)$ be the interior of $\cH(S^2)$ in the $C^2$ topology.
 One has to prove that $\cF^1(S^2)$ is empty.
 In~\cite{CP2} it is proven that if $g\in\cF^1(S^2)$ then the
 closure $\oPg$ of the set of periodic orbits for $g$ is a uniformly
 hyperbolic set. Moreover, it contains a non-trivial basic set
 $\La$. One can also assume that the geodesic flow for $g$
 is Kupka-Smale because in~\cite{CP2} it is proven that
 the Kupka-Smale metrics are  $C^2$ dense. 

 If the metric $g$ has positive curvature, George Birkhoff~\cite{Birk2}
 shows that there exists a simple closed geodesic $\ga$ such that
 the geodesic flow on $T^1S^2\setminus T\ga$ admits a global transversal 
 section homeomorphic to an open annulus. The Poincar\'e's return 
 map $f$ to the section preserves a finite area form and the return
 time is uniformly bounded away form zero and infinity.
 In section~\S\ref{proofs} we show how to use that $\oPg$
 is infinite, uniformly hyperbolic and that  $f$ is Kupka-Smale
 and area preserving to overcome the closing lemma problem and 
 arrive to a contradiction.

 When the curvature of $g$ is not strictly positive,
 the return map to the Birkhoff's section may not
 be defined everywhere. Instead of using the Birkhoff
 section we chose to lift the geodesic flow to a
 Reeb flow of a (tight) contact form on $S^3$.
 If the metric is Kupka-Smale, Hofer, Wysocki, Zehnder~\cite{HWZ1}
 theory implies that either 
 \begin{itemize}
 \item there is a periodic orbit $\ga$ such that the 
 Reeb flow on $S^3\setminus\ga$ admits a global transversal section
 homeomorphic to a disc with Poincar\'e's return map preserving
 a finite area and uniformly bounded return times. This
 is called the dynamically convex case and the proof in this case
 also holds for the Birkhoff's section in the positive curvature case.

\item[or]

 \item there is a finite set of periodic orbits $\ga_1,\ldots,\ga_N$
 and a global system of transversal sections
 on $S^3\setminus\cup_{i=1}^N\ga_i$.
 The return map preserves an area form of finite total area
 but it is not defined everywhere. The return time is bounded away from zero
 but it becomes infinite only at the intersection of the stable
 manifolds of a fixed finite set of periodic orbits.
 This case is called non dynamically convex.
 \end{itemize}
 
 In section~\S\ref{lift} we show how to lift the geodesic flow
 to the Reeb flow of a tight contact form on $S^3$.
 In section~\S\ref{s3} we summarize the Hofer, Wysocki, Zehnder
 theory on global sections for  generic Reeb flows of
 tight contact forms on $S^3$ and prove the lemmas that we need
 for the non dynamically convex case. In section~\S\ref{proofs}
 we prove theorem~\refteo.

 G. Contreras wants to thank the hospitality of 
 Universidade Federal de Minas Gerais where part of this research
 was done.

 % In~\cite{CP2} the Kupka-Smale theorem for geodesic flows is
 % proven: there is a  residual set $\cR$ of riemannian metrics with 
 % the $C^2$ topology such that for all $g\in\cR$ the closed geodesics
 % are either hyperbolic or generic elliptic and when the invariant 
 % manifolds of hyperbolic geodesics intersect, they intersect 
 % transversally.

 %By now, the main difficulty to generalize Newhouse theorem to riemannian
 %metrics is the lack of a $C^2$ closing lemma for geodesic flows.
 %The proof of Herman's theorem has two ingredients. The analogous to 
 %Ma\~ne's theory on dominated splittings for geodesic flows on 
 %surfaces, developed in~\cite{CP2}, that we shall briefly recall,
 %and the use of Brouwer's translation theorem applied to the return map
 %to the Birkhoff surface of section in order to bypass the closing lemma.

 \section{Lifting of the geodesic flow.}\label{lift}

 A contact manifold is a pair $(N,\la)$, where $N$ is an 
 odd dimensional smooth manifold and $\la$ is a contact 1-form,
 i.e. $\la\wedge d\la^n$ is a volume form on $N$, where $\dim N=2n+1$.
 The Reeb vector field $Y$ of $(N,\la)$ is the unique vector field  
 determined by 
 $$
 i_Y\,d\la \equiv 0\qquad \text{and} \qquad \la(Y)\equiv 1.
 $$

 Write $S^3=\{z\in\co^2\;|\;|z|=1\,\}$, $z=(z_1,z_2)=(q_1+ip_1,q_2+ip_2)$. 
 The standard contact form on $S^3$ is 
 $$
 \la_0:= \tfrac 12 {\textstyle\sum\limits_{\a=1}^2}
  \big[p_\a\,d q_\a-q_\a\, dp_\a\big]\Big|_{S^3}
  =\tfrac 12\, \big[ p\,dq- q\,dp\big]\vert_{S^3}.
 $$
 
 \begin{proposition}\label{conformal}
 For any riemannian metric $g$ on $S^2$, the double cover
 of the geodesic flow of $(S^2,g)$ is conjugate to 
 the Reeb flow of a positive multiple $f(z)\, \la_0$,
 $f(z)>0$ of the standard contact form $\la_0$ on $S^3$.
 \end{proposition}

 \begin{proof}

 We need to see the geodesic flow on $S^2$ as a Reeb vector field.
  For background material see~\cite{Klingenberg} or \cite{Patb}.   
 Let $\pi:TS^2\to S^2$ be the projection $\pi(x,v)=x$.
 The tangent space to $TS^2$ can be decomposed 
 $T_{\theta}TS^2=H(\theta)\oplus V(\theta)$ 
 where $V(\theta)=\ker d_\theta\pi$
 and $H(\theta)$ is the kernel of the connection map.
 The maps $d_\theta\pi: H(\theta)\to T_{\pi(\theta)}S^2$, 
 $\nabla_\theta:V(\theta)\to T_{\pi(\theta)}S^2$ are
 linear isomorphisms and they induce the Sasaki's  riemannian metric on 
 $TS^2$ by
 $$
 \langle\!\langle \zeta_1,\zeta_2\rangle\!\rangle_\theta:=
 \langle h_1,h_2\rangle_{\pi(\theta)}
 +
 \langle v_1,v_2\rangle_{\pi(\theta)},
 $$
 where $\zeta_i=(h_i,v_i)\in H\oplus V$.
 The Liouville's 1-form $\Theta$ on $T S^2$ is defined
  $\Theta_{(x,v)}\zeta:=\langle v,d\pi(\zeta)\rangle_x$,
 where $\pi:T^1S^2\to S^2$ is the projection $\pi(x,v)=x$.
 Its differential $\om=d\Theta$ is a symplectic form
 on $TS^2$ and is computed as 
 $$
 \om(\zeta_1,\zeta_2)=\langle v_1,h_2\rangle_x-\langle v_2,h_1\rangle_x,
 $$
 where $\zeta_i=(h_i,v_i)\in H\oplus V$. In particular $\om\wedge \om$ 
 is a volume form in $TS^2$. The vector field $Z(x,v)=(0,v)$ is the
 unit normal vector to the unit tangent bundle $T^1 S^2$ under 
 the Sasaki metric. Its contraction $i_Z\om$ is the Liouville's 
 form $\Theta$. Then $\Theta\wedge d\Theta=\tfrac 12 \,i_Z(\om\wedge\om)$
 is a volume form on $T^1S^2$ and $(T^1S^2,\Theta)$ is a
 contact manifold.

 In the decomposition $T(TS^2)=H\oplus V$,
 the geodesic vector field  is written as $X(x,v)=(v,0)$. 
 Then on the unit tangent bundle
 $T^1S^2$ we have that 
 \begin{gather*}
 \Theta(X(x,v))=\langle v, v\rangle_x\equiv 1,\\
 (i_X \om)_\theta(\zeta) =\langle \theta,d\pi(\zeta)\rangle_{\pi(\theta)}
 = \langle\!\langle Z(\theta),\zeta\rangle\!\rangle_\theta
 =0.
 \end{gather*}
 Then the geodesic vector field is the Reeb vector field of the
 Liouville's 1-form.

 % We need to see the geodesic flow on $S^2$ as a Reeb vector field.
 % The geodesic flow is the Euler-Lagrange flow of the kinetic energy
 % lagrangian $L:TS^2\to \re$, $L(x,v)=\frac 12 \,\lV v\rV_x^2$.
 % The Liouville's 1-form $\Theta$ on cotangent bundle $T^*S^2$ 
 % is defined by $\Theta_{(x,p)}(\zeta)= p(d\pi(\zeta))$, where 
 % $\pi:T^*S^2\to S^2$ is the projection $\pi(x,p)=x$.
 % Its differential $\om=d\la$ is the canonical symplectic 
 % form on $T^*S^2$.
 % The Legendre transform $\cL(x,v)=\langle v,\cdot\rangle_x$
 % is a conjugacy between the geodesic flow and the hamiltonian flow
 % of the hamiltonian $H:T^*S^2\to \re$, $H(x,p)=\tfrac 12\,\lV p\rV_x^2$,
 % where $\lV \cdot\rV_x$ is the induced metric on 
 % the cotangent bundle $T^*S^2$ by the riemannian metric on $S^2$.
 % The hamiltonian vector field $X$ for $H$ is defined by
 % $i_X\,\om=-dH$. 

 The unit sphere bundle $T^1S^2$ of $S^2$ is diffeomorphic to
 the special orthogonal group $SO(3)$
 by  identifying  the orthogonal matrix with columns
 $[x,v,x\times v]$ with the unit tangent vector $(x,v)$.
 The map $f:SO(3)\to S^2$, $f([x,v,x\times v])=x$ corresponds to
 the projection $\pi:T^1S^2\to S^2$.
 The fibers of the map $g([x,v,x\times v])=x\times v$
 are the unit tangent vectors to the oriented great circles of $S^2$
 (with axis $x\times v$).

 We show that the
  double cover of $SO(3)$ is the 3-sphere $S^3$.
 We identify  $S^3$  as the unit norm quaternions and
 $S^2$ as the unit norm quaternions $q\in\Q$ with zero real part.
 $$
 \Q=\{x_0+x_1i+x_2j+x_3k\;|\;x_\a\in\re,\;i^2=j^2=k^2=-1,\;
 ij=k,\; jk=i,\; ki=j\,\}.
 $$ 
 The covering map is $R:S^3\to SO(3)$, 
 $R_q(x):=q^{-1}\,x\,q=\ov{q}\,x\,q$.
 Indeed, if $|x|=1$, then $|R_q(x)|=1$ so that $R_q$ is orthogonal.
 Since $R_q(1) = 1$, $\rho(q)$ preserves the above embedding 
 $S^2\subset S^3$. Since $R_1= id$ and $q\mapsto R_q$ is
 continuous, we have that $R_q$ preserves orientation.
 The group $SO(3)$ is the set of orientation preserving
 rotations. It is generated by the sets $U_i$, $U_j$, $U_k$ 
 of rotations around the  three coordinate axis. 
 In order to prove that the map $R_q$ is surjective it is enough 
 to prove that $U_i$, $U_j$, $U_k$ are in the image of $R$.
 We only prove that $U_i\subset R(S^3)$.
 The rotation axis of $R_q$, $q=a+y$, $a\in\re$, $\Re(y)=0$
 is the direction of  $y$ in $\re^3=\Re^{-1}\{0\}\subset\Q$,
 because in this case $q\,y = y\, q$ and then $R_q(y)=\ov{q}\, y\, q=q$.
 If $q=a+ib=e^{i\theta}$ and  $x=zj$, with $z\in \co$,
 then 
 $ R_q(zj)=e^{-i\theta}zj\,e^{i\theta}
          =e^{-2i\theta} z j$. 
 So that $R_q$ is the rotation of angle $2\theta$ with axis $i$.

  Define the map $F:S^3\to SO(3)=T^1S^2$
 by $F(q)=[R_q(j),R_q(k),R_q(i)]=[x,v,x\times v]$.
 The lift under $F$ of the tangents to the oriented 
 great circles of $S^2$ are those $q\in S^3$ for which $R_q(i)$
 takes a fixed value.
 We show that $F^*\Theta_0=2\,\la_0$, where $\Theta_0$ is the 
 Liouville's form for the metric $g_0$ with  curvature  $+1$.

  Let $\langle\cdot,\cdot\rangle$ be the euclidean inner product on
 $\re^4\approx\Q$. Observe that for $q_1,\,q_2\in\Q$ we have that
 $$
 \langle q_1,q_2\rangle=\Re(q_1\,  \ov{q_2})=\Re(\ov{q_1}\, q_2)
 =\Re(q_2\,\ov{q_1}).
 $$
 For $\zeta\in T_qS^2$, $q=q_1+ p_1\, i+q_2\,j+p_2\,k$,
 $\zeta=x+y\,i+z\,j+w\,k$, the standard contact form $\la_0$ is written
 as
 \begin{align*}
 \la_0(q)\cdot\zeta &= p_1\, x-q_1\,y+p_2\,z - q_2\,w \\
     &=\langle p_1-q_1 i+p_2 j-q_2 k,\zeta\rangle 
     =\langle iq,\zeta\rangle \\
     &=\Re[-\ov{\zeta}\,i\,q].
 \end{align*}
 Since $\pi\scirc F(q)=\big(R_q(j),R_q(k)\big)=(x,v)$, 
 the first component of $d\pi\scirc d_q F(\zeta)$ is
 \begin{align*}
 d_qF(\zeta) =\big(d_q R(j)\cdot \zeta, *\big)
             =(\ov{\zeta}\,j\,q+\ov{q}\,j\,\zeta,*).
 \end{align*}
 The pull-back $F^*\Theta_0$ is given by
 \begin{align*}
(F^*\Theta_0)_q\cdot\zeta
  &=\Theta_0(\ov{q})\cdot d_qF(\zeta)
   =\langle R_q(k),\,\ov{\zeta}\,j\, q+\ov{q}\,j\,\zeta\rangle
   \\
   &=\langle \ov{q}\,k\,q,\,\ov{\zeta}\,j\,q+\ov{q}\,j\,\zeta\rangle
   =\Re\left[ (\ov{\zeta}\,j\,q)\,\ov{\ov{q}\,k\,q}
       +\ov{\ov{q}\,k\,q}\,(\ov{q}\,j\,\zeta)\right]
   \\
   &=\Re\left[\ov{\zeta}\,j\,q\,\ov{q}\,(-k)\,q
       -\ov{q}\,k\,q\;\ov{q}\,j\,\zeta\right]
   \\
   &=\Re\big[-\ov{\zeta}\,i\,q\big] +\Re\big[\ov{q}\,i\,\zeta\big].
 \end{align*}
 Since $\Re[\;\ov{q}\,i\,\zeta\;]=\Re[\;\ov{\ov{q}\,i\,\zeta}\;]
 =\Re[-\ov{\zeta}\,i\,q\,]$, 
 $$
 \big(F^*\Theta_0)_q\cdot\zeta= 2\,\Re\big[-\ov{\zeta} \,i\, q\,\big]
 = 2\,\la_0.
 $$

 If $g$ is another riemannian metric on $S^2$, using isothermal coordinates
 (e.g.~\cite{Chern1}) one can show that there exists a diffeomorphism
 $h:S^2\to S^2$ and a smooth positive function $f:S^2\to\re^+$
 such that $h^* g=f g_0$. Since $h$ is an isometry between
 $(S^2,fg_0)$ and $(S^2,g)$, the map $dh:TS^2\to TS^2$ conjugates
 their geodesic flows. Using $h$ we can assume that $g= f\, g_0$.
 Let $H:(T^1S^2,g_0)\to (T^1S^2,g)$ be the map 
 $(x,v)\mapsto \big(x,v/{\scriptstyle\sqrt{f(x)}}\big)$ and
 let $\Theta$ be the Liouville's form for the metric $g=f\,g_0$.
 Then
 \begin{align*}
 H^*\Theta_{(x,v)}\cdot\zeta
  &= \Theta_{H(x,v)}\cdot d_{(x,v)}H(\zeta) \\
  &=\Theta\big(\tfrac{v}{\sqrt{f}}\big)\cdot d\pi(\zeta)\\
  &=g_x\big(\tfrac{v}{\sqrt{f}},\,d\pi\cdot\zeta\big)
   =f(x)\;g_0\big(\tfrac{v}{\sqrt{f}},\,d\pi\cdot\zeta\big)\\
  &=\sqrt{f(x)}\,g_0(v,d\pi\cdot\zeta)
   =\sqrt{f(x)}\;\Theta_0{(x,v)}\cdot\zeta.
 \end{align*}
 So that $H^*\Theta=\sqrt{f\scirc\pi}\;\,\Theta_0$.
 Then $F^*H^*\Theta=2\,\sqrt{f\scirc\pi\scirc F}\,\,\la_0$.
 If $X$ is the vector field of the geodesic flow for $(S^2,g=f\,g_0)$,
 since $X$ is the Reeb flow of the Liouville's form for $g$,
 then $(F\scirc H)^{-1}(X)$ is the Reeb vector field for 
 $2\,\sqrt{f\scirc\pi\scirc F}\,\,\la_0$.
 \end{proof}

 \section{Generic tight contact flows on $S^3$}\label{s3}

 In this section we summarize the theory of contact flows on $S^3$ by Hofer, 
 Wysocki and Zehnder that we shall need. A description of 
 the theory is given in~\cite{HWZ2} and the proofs can be found 
 in~\cite{HWZ1}.
 
 Write $S^3:=\{\,z\in\co^2\,|\,|z|=1\,\}$, 
 $z=(z_1,z_2)=(q_1+i p_1, q_2+i p_2)$  with $z_j\in\co$ and $p_j,\,q_j\in\re$.
 Let 
 $$
 \la_0:=\tfrac 12\,\sum_{j=1}^{2}\big[ q_j\,dp_j-p_j\, dq_j\big]\Big\vert_{S^3}
 $$
 be the standard contact form in $S^3$. 
 In the following we shall only consider contact forms
 $$
 \la= f\cdot \la_0,
 $$
 where $f:S^3\to]0,+\infty[$ that we shall call tight contact forms.
 Any such $\la$ is indeed a contact form in $S^3$ because
 $$
  \la\wedge d\la = f^2\cdot (\la_0\wedge d\la_0)
 $$
 is a volume form if $f^{-1}\{0\}\ne\emptyset$.
 The {\it Reeb vector field } of $\la$ is the vector field $X$ on $S^3$ 
 defined by
 $$
 i_X\, d\la\equiv 0
 \quad\text{and}\quad \la(X)\equiv 1.
 $$
 Its flow $\vr_t$ is called the {\it Reeb flow} of $\la$. 
 The {\it contact structure} of $\la$ is the distribution of 
 linear 2-dimensional subspaces
 $$
 \xi =\ker \la.
 $$
 The contact form $\la$ and its contact structure $\xi$ are invariant
 under the derivative of the Reeb flow because
 $$
 L_X \la = d (i_X \la) + i_X d\la \equiv 0.
 $$
 The derivative $d\la$ is a symplectic (i.e area-) form on any 2-dimensional
 subspace which is transversal to the vector field $X$, in particular
 on the contact structure $\xi$.
 
 \subsection{Closed orbits of the Reeb flow}\quad
 \medskip
 
 Let $x:[0,T]\to S^3$ be a periodic orbit for $\vr_t$ with  period $T$.
 Then $T$ is a multiple of the minimal period of $x$.
 We shall say that $(x,T)$  is {\it non-degenerate} if
 the number $+1$ is not an eigenvalue of its linearized Poincar\'e map
 $d\vr_T\vert_\xi$, restricted to $\xi=\ker \la$. 
 Since $d\vr_T\vert_\xi$
 preserves the area form $d\la\vert_\xi$, its eigenvalues have the form
  $\{\mu_1,\mu_2\}=\{\mu,\ov{\mu},\mu^{-1},(\ov{\mu})^{-1}\,\}$.
 Hence, if $(x,T)$ is non-degenerate, either $|\mu|=1$ or
  $\mu\in\re\setminus\{0\}$. We say that $(x,T)$ is {\it hyperbolic}
 if $\mu\in\re$, $\mu\ne \pm 1$ and that $(x,T)$ is {\it elliptic }
 if $\mu\in\co\setminus\re$, $|\mu|=1$.
 We shall also distinguish between (+)-hyperbolic orbits, when 
 $\mu, \mu^{-1}>0$;
 and ($-$)\,-hyperbolic orbits, when $\mu, \mu^{-1}<0$.
 
 We describe now one characterization of the Conley-Zehnder index 
 of a non-degenerate periodic orbit $(x,T)$. Assume that $(x,T)$ 
 is non-degenerate. Since $\pi_1(S^3)=0$,  
 \linebreak
 $x:\re/{\scriptstyle T\Z}\to S^3$ is contractible.
 Choose a map $v_D:D\to S^3$,  $D:=\{\,z\in\co\,|\,|z|\le 1\,\}$
 such that $v_D(e^{2\pi i s})=x(Ts)$.
 Then the contact bundle $v_D^*$ is trivial. Choose a
 trivialization $v_D^*\to D\times \re^2\approx D \times \co$.
 Write the derivative of the Reeb vector field on this trivialization
 $$
 \Phi(s):=d\vr_{sT}\vert_\xi\in \Sp(1)
 =\{\,A\in\re^{2\times 2}\,\vert\,\det A=+1\,\}\;,
 \qquad 0\le s\le 1.
 $$
 This arc of symplectic matrices starts at the identity $\Phi(0)=I$ 
 and ends at the linearized Poincar\'e map $\Phi(1)=d\vr_t\vert_\xi$.
 Let $z\in\co\setminus\{0\}$ and $z(s):=\Phi(s)\,z$.
 Choose a continuous argument $\theta(s)$ for 
 $$
 e^{2\pi i\theta(s)}=\frac{z(s)}{|z(s)|},
 \qquad 0\le s\le 1.
 $$
 Define the winding number of $\Phi(s)z$ by
 $$
 \Delta(z):=\theta(1)-\theta(0)\in\re,
 $$
 and the winding interval of the arc $\phi$ by
 $$
 I(\Phi):=\big\{\,\Delta(z)\,\big\vert\;z\in\co\setminus\{0\}\,\big\}.
 $$
 \begin{lemma}
 \quad
 $\length\big[I(\Phi)\big]\le\tfrac 12$.
 \end{lemma}
 \begin{proof}
 Let $z(s):=\Phi(s) z$, $w(s):=\Phi(s)\,w$ and $u(s):=z(s)\, \ov{w(s)}\in\co$.
 Observe that $\Delta(u)=\Delta(z)-\Delta(w)$. Assume that
 $|\De(u)|\ge \tfrac 12$. Then there exists $0<s_0<1$ such that 
 $z(s_0)=\tau\,w(s_0)$  for some $\tau\in\re\setminus\{0\}$.
 Then $z(s)=\tau\,w(s)$ for all $s\in[0,1]$. Hence
 $\De(z)=\De(w)$, which contradicts $|\De(z)-\De(w)|\ge\tfrac 12$.
 \end{proof}
 
 Then the winding interval either lies between two consecutive integers 
 or contains an integer.
 Define the {\it Conley-Zehnder index } of the non-degenerate periodic 
 orbit $(x,T)$ by
 $$
 \mu(x,T):=\mu(\Phi):=\left.\begin{cases}
 2k+1 & \text{if }\; I(\Phi)\in]k,k+1[\; \\
 2k & \text{if }\;k\in I(\Phi)
 \end{cases}
 \right\},
 \quad k\in\Z. 
 $$
 The integer $\mu(\Phi)$ depends only on the homotopy type of the chosen 
 disc map $v_D:D\to S^3$. Since $\pi_2(S^3)=0$, 
 the index\footnote{i.e. two such discs
 $v_i:D_i\to S^3$, $i=0,1$, can be joined to form a sphere 
 $u:S^2=D_0\cup D_1\to S^3$. 
 Since $\pi_2(S^3)=0$, 
 $u(S^2)$ is the boundary of a 3-ball $B^3$ in $S^3$. 
 One can use the ball $B^3$ to construct a homotopy 
 $v_t$, between $v_0$ and $v_1$ with $v_t(e^{2\pi i s})=x(sT)$ for all 
 $s,t\in[0,1]$.} 
  $\mu(x,T)$ is well defined.
 
 %\medskip 
 
 Observe that the winding number $\De(z_0)$ is an integer 
 if and only if $\Phi(1) z_0=\mu z_0$ for some $\mu>0$.
 Hence
 \begin{alignat*}{3}
 &\mu(x,T) \text{ is even }\quad
 &&\Longleftrightarrow 
 &&\quad\; (x,T) \text{ is (+)-hyperbolic.}
 \\
 &\mu(x,T) \text{ is odd }\quad
 &&\Longleftrightarrow 
 &&\quad\;(x,T) \text{ is elliptic or ($-$)-hyperbolic.}
 \end{alignat*}
 
 \bigskip
 
 We also define the self-linking of a closed orbit $(x,T)$.Take a disc map 
 $v_D:D\to S^3$ such that $v_D(e^{2\pi i s})=x(sT)$, $0\le s\le 1$. Choose
 a nowhere vanishing section $Z$ of the pull-back bundle $v_D^*\xi$.
 Then $Z$ is nowhere tangent to $x$. Pushing $x$ slightly in the direction of
 $Z$ we obtain a loop $x'$ which is transversal to $\xi$ and disjoint from $x$.
 The two loops $x$, $x'$ have a natural orientation induced by the 
 orientation of $\partial D$. Define the 
 {\it self-linking number of $(x,T)$} 
 $$
 sl(x,T):= I(x',v_D)\in \Z,
 $$
 as the oriented intersection number of $x'$ with $v_D$.
 This number does not depend on the choices of $v_D$ or $Z$.
 Indeed, using a trivialization of $v_D^*\xi$, the section $Z$ 
 is a function $Z:D\to\co\setminus\{0\}$. Since 
 $\pi_2(\co\setminus\{0\})=0$, $Z$ is homotopic to a constant
 function $Z:D\to\{1\}$. Then the loop $x'$ is isotopic to the loop 
 $\{\ov x\}$ obtained from $Z$ and thus it has the same 
 intersection number as $x$. Similarly, $sl(x,T)$ only depends
 on the homotopy type of $v_D$, but since $\pi_2(S^3)=0$, 
 there is only one homotopy type.
 
 \subsection{Finite energy surfaces}\quad
 
 \medskip
 
 An {\it almost complex structure compatible with $\la$} on $S^3$
 is a linear bundle map $J:\xi\to\xi$ such that $J^2=-Id$ and such that
 the quadratic form 
 $$
 \xi\times\xi\ni(h,k)\longmapsto d\la(h,Jk)
 $$
 is positive definite on each fiber of the product bundle $\xi\times\xi$.
  The set of smooth 
 almost complex structures compatible with $\la$ is always non-empty
 and contractible. 
 We extend $J$ to an almost complex structure $\tJ$ on 
 the product $\re\times S^3$ by setting
 \begin{align*}
 \tJ(a,bX+h)=\big(\!-b\;,\;aX+Jh \big), \quad\text{ for }\;\; a,b\in\re,\;
 h\in\xi.
 \end{align*}
 
 Let $\Si:=\co\cup\{\infty\}$ be the Riemann sphere and let $\Ga\subset \Sigma$
 be a finite set (of ``punctures''). 
 Let $j:T\Si\to T\Si$, $j(z)=iz$, be the complex structure on $\Si$.
 
 \begin{definition}[{cf. \cite[appendix]{HWZ2}}]\quad
 
 \noindent
 A (spherical) {\bf finite energy surface} is a map $\tu=(a,u):\SiGa\to\reS$
 such that
 \begin{gather}
 \tu:\SiGa\to\re\times S^3 \text{ is proper and non-constant,} 
 \label{pnc} \\
 d\tu\scirc j=\tJ\scirc d\tu, 
 \label{psh}\\
 \int_\SiGa u^*d\la<+\infty.
 \label{dle}
 \end{gather}
 A {\bf finite energy plane} is a spherical finite energy surface
 with only one puncture $\Ga=\{\infty\}$.
 \end{definition}
 
 The integrand in~\eqref{dle} is always non-negative, and it is positive
 at the points where the projected map $u:\SiGa\to S^3$ is transversal to
 the Reeb vector field.
 
 The set of punctures $\Ga$ must be non-empty because
 if $\Ga=\nulo$ then $\tu$ is constant~\cite[lemma~3.4]{HWZ2}.
 The behaviour of a finite energy surface near a puncture
 is classified by the following
 \begin{Lemma}[\cite{HWZ2}]
 Let $\tu:\SiGa\to\reS$ be a finite energy surface and $\ga\in\Ga$ 
 a puncture. Then one of the following cases hold, where $\tu=(a,u)$:
 \begin{itemize}
 \item positive puncture: \;\; $\lim\limits_{z\to\ga}a(z)=+\infty$;
 \item negative puncture: \,\, $\lim\limits_{z\to\ga}a(z)=-\infty$;
 \item removable puncture:  $\lim\limits_{z\to\ga}a(z)$ {\rm exists in }$\re$.
 \end{itemize}
 \end{Lemma}
 
 In the case of a removable puncture $\ga$ there exists a 
 neighbourhood $\cU(\ga)$ of $\ga$ such that $\tu$ is bounded 
 on $\cU(\Ga)\setminus\{\ga\}$. By Gromov's removable singularity 
 theorem~\cite{Gro30}, $\tu$ can be extended smoothly over the 
 puncture $\ga$. Hence we only consider positive and negative punctures:
 $$
 \Ga =\Ga^+ \cup \Ga^-.
 $$ 
 
 A special example of a finite energy surface is the following.
 Let $(x,T)$ be a periodic orbit of the Reeb vector field. 
 Then the map $\tu:\co\setminus\{0\}\to\reS$,
 $$
 \tu\big(e^{2\pi(s+it)}\big)=\big(\,Ts,x(Tt)\big)\in \re\times S^3.
 $$
 is an embedded finite energy surface in $\reS$. 
 Its  $d\la$-energy
 vanishes:
 $$
 \int_{\co\setminus\{0\}}u^* d\la=0.
 $$
 Its image $F=\tu(\co\setminus\{0\})$ is  fixed under the $\re$-action
 on $\re\times S^3$.

 At a puncture $\ga$, a finite energy surface converges to a
 periodic orbit $(x_\ga,T_\ga)$ of the Reeb vector field
 which is called the {\it asymptotic limit} of
 the surface at the puncture. 
 Its period $T_\ga$ can be a multiple of its minimal period $\tau_\ga$.
 When $T_\ga=\tau_\ga$ we say that $(x_\ga,T_\ga)$ is simply covered by $\tu$.
 Now we describe better the local form  of $\tu=(a,u):\SiGa\to\re\times S^3$ 
 near a puncture.

  Assume that the asymptotic limit $(x,T)$ associated with the puncture
 is a  non-degenerate periodic orbit.  
  Fix holomorphic polar coordinates $\si(s,t)$
 in a punctured neighbourhood of the puncture $\ga\in \Ga$, 
 we have the asymptotic behaviour
 $$
 \lim_{s\to\infty}u\scirc\si(s,t)=x(mt) \quad\text{ in } C^\infty(S^1),
 $$
 where $\ga=\lim_{s\to\infty}\si(s,t)$ and $(x,T)$ is a periodic 
 orbit of the Reeb vector field with (perhaps non-minimal) period $T=|m|$.
 When $\ga\in\Ga^+$,  $m=T$
 and when $\ga\in\Ga^-$, $m=-T$. 
 
 One can give coordinates to a neighbourhood of 
 the asymptotic limiting periodic
 orbit $P$ of the form $S^1\times\re^2$, $S^1=\re/\Z$,
  with $P=S^1\times\{0\}$ and
 $\{0\}\times\re^2$ tangent to the constant structure $\xi$.
 In such coordinates 
 \begin{itemize}
 \item[either]
 \item there exists $c\in\re$ such that
 $$
 \tu\scirc\si(s,t)=(ms+c,mt), \qquad \text{ for }(s,t)\in S^1\times\re^2,
 $$ 
 and its whole projection $u(\SiGa)=P$ is only the asymptotic limit, 
 as in the example above;
 
 \item[or,] 
 
 \item the projected map $u$ has the form
 $$
 u\scirc\si(s,t)\approx \big(mt, e^{\mu s} e(t)\big), 
 \qquad\text{as }s\to\infty;
 $$
 modulo lower order terms in $s$. 
 In  a positive puncture $s\to+\infty$ and $\mu<0$, in
 a negative puncture $s\to-\infty$ and $\mu>0$.
 Therefore near such a puncture there is a precise directional convergence 
 of $u$ to the periodic orbit given by the periodic
 non-vanishing vector field 
 \begin{equation}\label{et}
 e(t)=\lim_{s\to\infty}
 \frac{\partial_s u}{\left|\partial_s u\right|}
 \in\xi_{x(mt)}\setminus\{0\},
 \qquad t\in S^1=\re/\Z
 \end{equation}
 in the contact structure $\xi$ along $P$.
 \end{itemize} 

 In the later case above, by Stokes theorem
 $$
 \sum_{\ga\in\Ga^+}T_\ga- \sum_{\ga\in\Ga^+}T_\ga
 =\int_\SiGa u^* d\la > 0,
 $$
 where $(x_\ga,T_\ga)$ is the asymptotic limit of $\tu$ at $\ga$.
 It follows that such  $\tu$ must have at least one positive puncture. 

 A neighbourhood of the puncture looks like $]0,+\infty[\times S^1$. 
 To this we add $\{+\infty\}\times S^1$. In this way one obtains 
 a circle compactification $\ov{\Si}$ of $\Si\setminus\{\ga\}$ and 
 $u:\SiGa\to S^3$ can be extended to a smooth map 
 $\ov{u}:\ov{\Si}\to S^3$ such that the boundary circles 
 parametrize the periodic orbits which are the asymptotic limits
 of the punctures. 
 On a positive puncture, $\ga\in\Ga^+$, the
 orientation of the limiting periodic orbit coincides with 
 the orientation of the boundary circle of $\ov{\Si}$ and 
 on a negative puncture the orientations are reversed.

 Define the {\it Conley-Zehnder index of $\tu:\SiGa\to\reS$} by
 $$
 \mu(\tu)=\sum_{\ga\in\Ga^+}\mu(x_\ga,T_\ga)
 -\sum_{\ga\in\Ga^-}\mu(x_\ga,T_\ga),
 $$
 where $(x_\ga,T_\ga)$ is the asymptotic limit of $u$ at $\ga$
 and $\mu(x,T)$ is the Conley-Zehnder index of a periodic orbit 
 $(x,T)$ defined above.
 
 Note that given a finite energy surface $\tuas$,
 and $c\in\re$,  the translated map
 $$
 \tu_c(z):=\big(a(z)+c,u(z)\big), z\in\SiGa,
 $$
 is also a finite energy surface because $\tJ$ 
 is invariant under translations in the first factor $\re$.
 Hence a finite energy surface always belong to
 a 1-parameter family of finite energy surfaces, invariant 
 under the $\re$-action on $\reS$.

  We define the {\it index} of an {\it embedded }finite energy surface
 $F=\tu(\SiGa)$ by
 $$
  \Ind(F):=\mu(F)-\chi(S^2)+\sharp F,
 $$
 where $\mu(F)=\mu(\tu)$, $\chi(S^2)=2$ the Euler characteristic 
 of $S^2$ and $\sharp F:=\sharp \Ga$. Observe that if $\tu$ is 
 an embedding $\mu(F)$ only depends on the image $F$ and not on 
 $\tu$. This index describes the dimension of the moduli space
 consisting of nearby embedded finite energy surfaces having 
 the same topological type.

 \subsection{Zoology of embedded finite energy surfaces with low indices}
 \quad
 
 \medskip
 
 In the following we shall only be interested in finite energy surfaces
 $$
 \tuas
 $$
  such that 
 \begin{itemize}
 \item its projection $\us$ is an embedding,
 \item have exactly one positive puncture, 
 \item its image $F=u(\SiGa)$ has indices
 $\mu(F)\in\{1,2,3\}$ and $\Ind(F)\in\{1,2\}$.
 \end{itemize}
  
 Assume that all the periodic orbits are non-degenerate.
 Let $(x,T)$ be an asymptotic limit of such $\tu$.
 Take a disc map $v_D:D\to S^3$·such that $v_D(e^{2\pi t})=x(tT)$ 
 and choose a nowhere vanishing section $Z$ of the pull-back 
 $v_D^*\xi$ of the contact bundle $\xi$.
 Then  there exists a nowhere vanishing function
 $f(t)\in\xi_{x(t)}$ along the orbit such that the linearized Reeb
 flow along the orbit has the representation
 $$
 d_{x(0)}\vr_t \cdot v=f(t)\cdot Z(x(t)),
 $$
 where the dot at the right denotes the complex multiplication 
 with respect to the complex structure $J$ on $\xi$. 
 By the definition of the Conley-Zehnder index we have the 
 following table
 
 \bigskip

 \centerline{
 \begin{tabular}[t]{|c|c|c|c|}
 \hline
 $\mu(x,T)$ & 
 $\begin{matrix}\text{change of argument of $f$} \\
 \text{along a full period}
 \end{matrix}$ &  
 \phantom{$\Big($} & eigenvalues of $d_{x(0)}\vr_T$ \\
 \hline
 \hline
 1 & $0<\De\arg f < 2\pi$ &\phantom{$\Big($}
 & elliptic or ($-$)-hyperbolic\\
 \hline
 2 & $\De\arg f =2\pi$ & 
 $\begin{matrix}\text{if $v$ is an eigenvector of} \\ 
  \text{the linearized Poincare map}
 \end{matrix}$
 & ($+$)-hyperbolic\\
 \hline
 2 & $\pi<\De\arg f < 3\pi$ & otherwise \phantom{$\Big($}
 & ($+$)-hyperbolic\\
 \hline
 3 & $0<\De\arg f < 2\pi$ &\phantom{$\Big($} 
 & elliptic or ($-$)-hyperbolic\\
 \hline
 \end{tabular}
 }
 
 \bigskip
 
 If the projection $u$ of $\tu$ is an embedding and $\ga\in\Ga$
 is a puncture there is a precise directional convergence in 
 $\xi$ towards the asymptotic periodic orbit, 
 described by the vector field $e(t)\in\xi_{x(t)}$ from~\eqref{et}.

 \parbox{7cm}{
 \includegraphics[scale=.35]{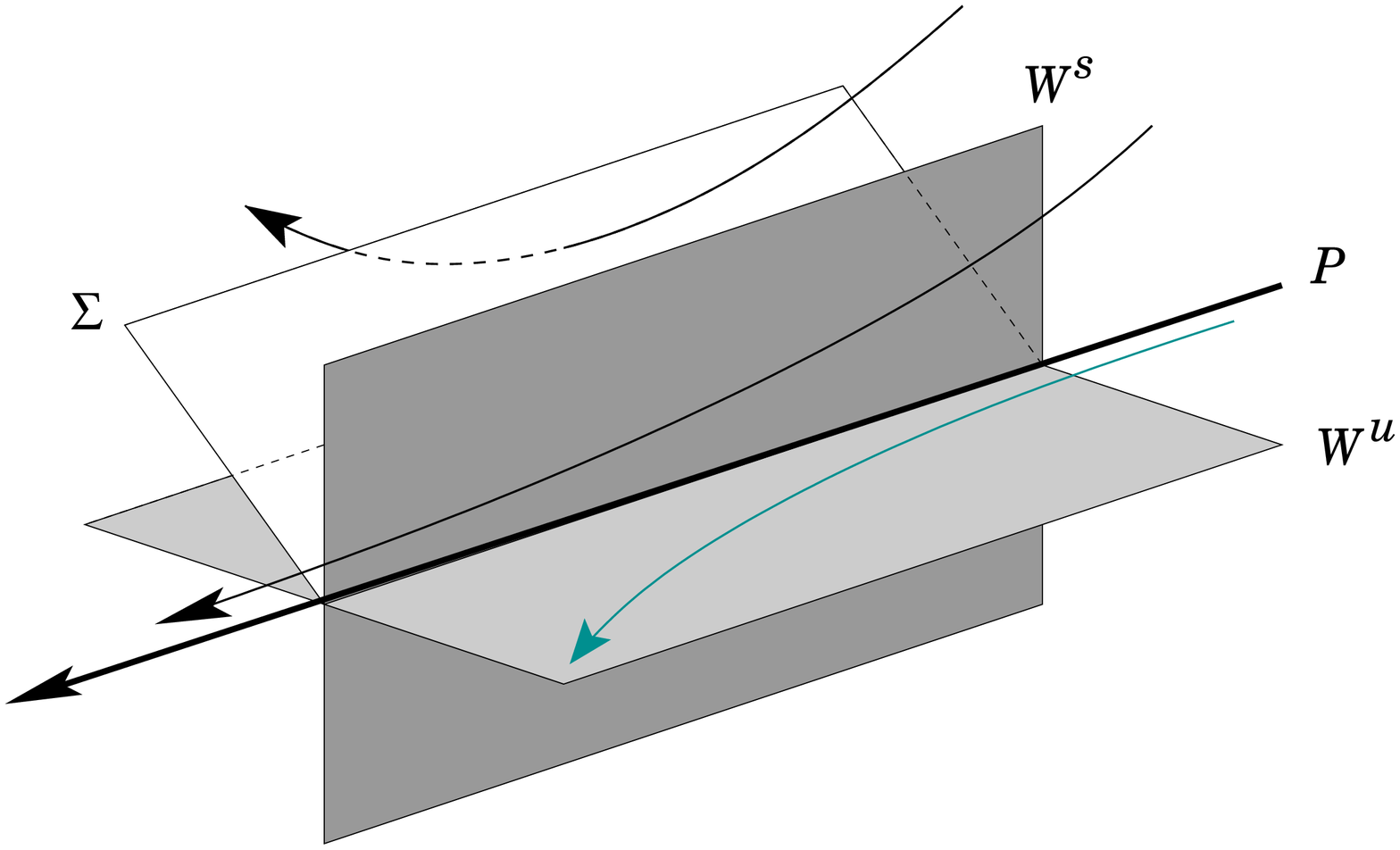}
 \refstepcounter{figure}\label{hyp}
 \\
 { \openup -2pt {\small {\sc Figure~\ref{hyp}.} 
 This figure shows the stable and unstable 
 manifolds of a binding periodic  orbit $P$ of Conley-Zehnder index 2
 and the approach of a finite energy surface $\Si$ having $P$ 
 as asymptotic limit.
 }}
 }
 \hskip 1cm
 \parbox{7cm}{
 \includegraphics[scale=.35]{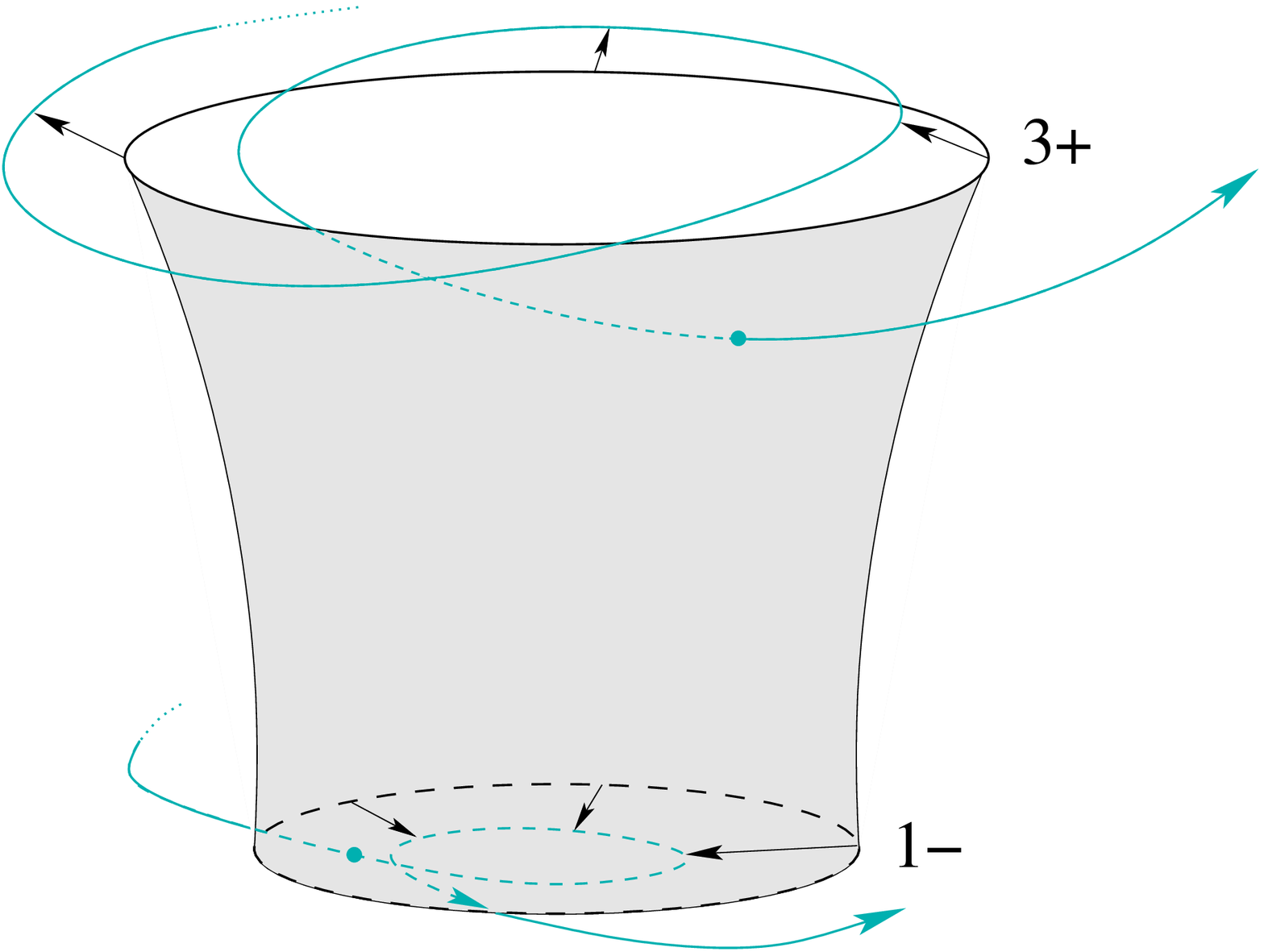}
 \refstepcounter{figure}\label{ell}
 \\
 { \openup -2pt{\small
 {\sc Figure~\thefigure.} The figure shows the winding of the flow 
 near periodic orbits of Conley-Zehnder indices 1 and 3, 
 relative to a finite energy cylinder. The Reeb vector field is 
 transversal to the finite energy surface. Two orbits intersecting
 the surface are shown.} }
 }
 
 \bigskip
 
 Consider for example $\ga\in\Ga^+$ and assume that the 
 asymptotic limit $P_\ga$ has  Conley-Zehnder index 2 or 3. Then 
 it can be proved (cf.~\cite{HWZ42,HWZ1}) that the vector $e$ 
 describing the approach towards the periodic orbit 
 has winding number at most 1 with respect to 
 the nowhere vanishing section $Z$.
 This means that the flow turns faster (more than $2\pi$) 
 around the periodic orbit than the approaching surface 
 if the Conley-Zehnder index is 3.
 
 If the Conley-Zehnder index is 2, the situation is a bit
 more subtle. Assume that the vector $e(t)$ describing the 
 directional convergence has winding number 1.
 Since $\mu(P_\ga)=2$, the  stable and unstable subspaces 
 of the periodic orbit
 have winding number 1 (they are generated by nowhere 
 vanishing sections of $\xi$ whose winding can be measured with
 respect to $Z$) and intersect a transversal section in a pair
 of transversal lines creating four quadrants. The trace of $e$ 
 on this transversal section will lie on one of this quadrants.
  See figures~\ref{hyp},~\ref{tranhyp},~\ref{loop} for more detail.
 
 If $\ga\in\Ga^-$ and if the asymptotic limit has Conley-Zehnder 
 index 1 or 2, then the winding number of the vector $e(t)$ 
 describing the asymptotic approach is at least 1. This time 
 the surface turns faster than the flow around the orbit if 
 the Conley-Zehnder index is 1.
 If the index is equal to 2, the periodic orbits is (+)-hyperbolic
 and the surface turns as fast as the flow, similar to the 
 behaviour near the positive puncture of index 2. 
  See figures~\ref{hyp},~\ref{ell},~\ref{tranhyp},~\ref{tranell}.
 
 \bigskip

 Now assume that $\tu$ satisfies the three items above. 
 Let $\mu^+$ be the Conley-Zehnder index of its (unique) 
 positive puncture. The negative punctures in $\Ga^-$ all have 
 $\mu$-indices~$\ge 1$ if $\Ga^-\ne \nulo$. Denote by $N_j$ 
 the number of negative punctures having $\mu$-index equal 
 to $j$, where $j\ge 1$. We have that $\Ind(F)\ge 1$.
 Hence
 $$
 \Ind(F)=\mu^+-\sum_{j=1}^\ell{j\,N_j}-2+(1+\sharp\Ga^-)\ge 1.
 $$
 This is equivalent to
 \begin{equation}\label{summu+}
 \sum_{j=1}^\ell(j-1)\,N_j\le\mu^+-2.
 \end{equation}
 We conclude for the index of the positive puncture that 
 $$
 2\le\mu^+.
 $$
 If $\mu^+=2$, it follows from~\eqref{summu+} that $N_j=0$ for all $j\ge 2$.
 If $\mu^+=3$, we conclude that $N_2\le 1$ and $N_j=0$ for all $j\ge 3$.
 In both cases there is no restriction on the number of negative 
 punctures with $\mu$-index equal to 1. In order to list the types 
 of  such surfaces $F$ we introduce the vectors
 $$
 (\mu^+,\mu_1^-,\ldots,\mu_N^-),
 $$
 where $N=\sharp \Ga^-$ is the number of negative punctures 
 of $F$, $\mu^+$ is the Conley-Zehnder index of the unique
 positive puncture and $\mu_j^-$ are the indices of the 
 negative punctures ordered so that $\mu_j^-\ge \mu_{j+1}^-$.
 Then $F$ must have one of the following types:
 \begin{alignat}{2}
 &(3,1_1,\ldots,1_N), \qquad &&\Ind(F)=2; 
 \notag \\
 &(3,2,1_1,\ldots,1_{N-1}), \qquad &&\Ind(F)=1;
 \label{tabla2}\\
 &(2,1_1,\ldots,1_N),\qquad &&\Ind(F)=1.
 \notag
 \end{alignat}
 The number $N$ of negative punctures can be zero.
 If this happens, the first and third cases represent finite energy planes.
 The second case, for $N=1$, represents a finite energy cylinder 
 connecting a periodic orbit of index 3 [elliptic or 
 ($-$)-hyperbolic] with a periodic orbit of index 2 [(+)-hyperbolic].
 
 We point out that the only index that can occur at a positive 
 and a negative puncture, although not simultaneously, is the 
 index equal to 2, which belongs to a (+)-hyperbolic orbit, 
 where the asymptotic approach of $F$ lies in a 
 quadrant between the stable an unstable manifolds.
 This will be important later on in the decomposition of families 
 of leaves at their ends  which always takes place along a 
 hyperbolic orbit of index 2.

 \subsection{Global system of transversal sections}
 
 \medskip
 
 \begin{definition}
 We say that a contact form $\la=f\la_0$ is {\bf non-degenerate}
 if all the periodic orbits $(x,T)$ of its Reeb vector field are
 non-degenerate, i.e. the eigenvalues of the linearized Poincare
 map for their minimal periods, restricted to $\xi=\ker \la$,
 are not roots of unity.
 \end{definition}
  We say that 
 
 \begin{Theorem}{\bf [Hofer, Wysocki, Zehnder~\cite{HWZ1,HWZ2}]}
 \label{HWZ}

 If  all the periodic orbits of the Reeb vector field $X$
 of the contact form $f\la_0$ are non-degenerate, then 
 there exists a non-empty set $\cP$ of finitely many periodic orbits
 distinguished periodic orbits of $X$ which  
 have self linking number $-1$ and Conley-Zehnder indices 
 in the set $\{1,2,3\}$ so that the complement
 $$
 S^3\setminus \cP
 $$
 is smoothly foliated into leaves which are embedded punctured 
 Riemann spheres, transversal to the Reeb vector field $X$ 
 and converging at the punctures to periodic orbits from $\cP$. 
 
 The leaves $F=u(\SiGa)$ are projections 
 of embedded finite energy surfaces with indices 
 \linebreak
 $\Ind(F)\in\{1,2\}$
 whose asymptotic limits (in $\cP$) are simply covered.
 Each leaf has precisely one positive puncture, but an arbitrary number
 of negative punctures.
 
 Moreover, there is at least one leaf which is the projection of a 
 finite energy plane and whose asymptotic limit $P$ has Conley-Zehnder 
 index $\mu(P)\in\{2,3\}$, {\rm (c.f.~\cite[prop. 7.1]{HWZ1}).} 
 
 \end{Theorem}
 
 \medskip
 
 Since the leaves $F=u(\SiGa)$ are transversal to the Reeb vector field,
 the 2-form $d\la\vert_F$ is a positive area form on $F$ and its total 
 area $\int_F d\la$ is finite.
 
 The periodic orbits in $\cP$ are called {\it binding
 periodic orbits} and the leaves $F$ with index $\Ind(F)=1$
 are called {\it rigid leaves}.
 
 The decomposition of $S^3$ in theorem~\ref{HWZ} comes from a foliation
 on $\reS$ by embedded finite energy surfaces. The binding periodic orbits
 are projections of surfaces $\tF=\tu(\SiGa)\subset\re\times P$, $P\in\cP$, 
 which are fixed under  the $\re$-action on $\reS$. 
 
 The rigid leaves
 are projections of surfaces $\tF=\tu(\SiGa)$, with $\Ind(\tF)=1$, 
 which belong to a 1-parameter family of surfaces having the same 
 asymptotic limits, namely, the orbit of $\tF$ under the $\re$ action
 on $S^3$. The projection of the family is an isolated 
 embedded punctured sphere.  There are finitely many binding orbits.
 
 The leaves $F=u(\SiGa)$ with $\Ind(F)=2$ belong to a 2-parameter
 family of finite energy surfaces all with the same asymptotic limits.
 One parameter is given by the orbit under the $\re$-action on $\reS$.
 After projecting by $\reS\to S^3$, it remains a 1-parameter family

 The foliation gives a decomposition
 $$
 S^3=\cS_0\cup\cS_1\cup\cS_2,
 $$
 where $\cS_0$ is the set of points in the binding periodic orbits,
 $\cS_1$ is the set of points in the finitely many rigid leaves 
 and $\cS_2$ is the set of points in the leaves with index $\Ind(F)=2$
 which occur in 1-dimensional families parametrized either by $S^1$ 
 or by an open interval $I=]0,1[$. In the case of an $I$-family $F_\tau$,
 the surfaces decompose at each of the ends where $\tau\to 0$ and $\tau\to 1$ 
 into two rigid surfaces $\{C^+,C^-\}$ along a hyperbolic binding 
 orbit $P$ whose Conley-Zehnder index is equal to $\mu(P)=2$.
 On $C^+$ the orbit $P$ is the asymptotic limit at its positive puncture
 and on $C^-$, $P$ is the asymptotic limit at a negative puncture.
 
 \begin{definition}
 We say that the contact form $f\la$ is {\bf dynamically convex}
 if it admits a foliation as in theorem~\ref{HWZ} which has
 an $S^1$-family of leaves $F_\tau$ with index $\Ind(F_\tau)=1$.
 \end{definition}

 \vskip 1cm
 
 \centerline{
 \parbox{5cm} {\includegraphics[scale=.38]{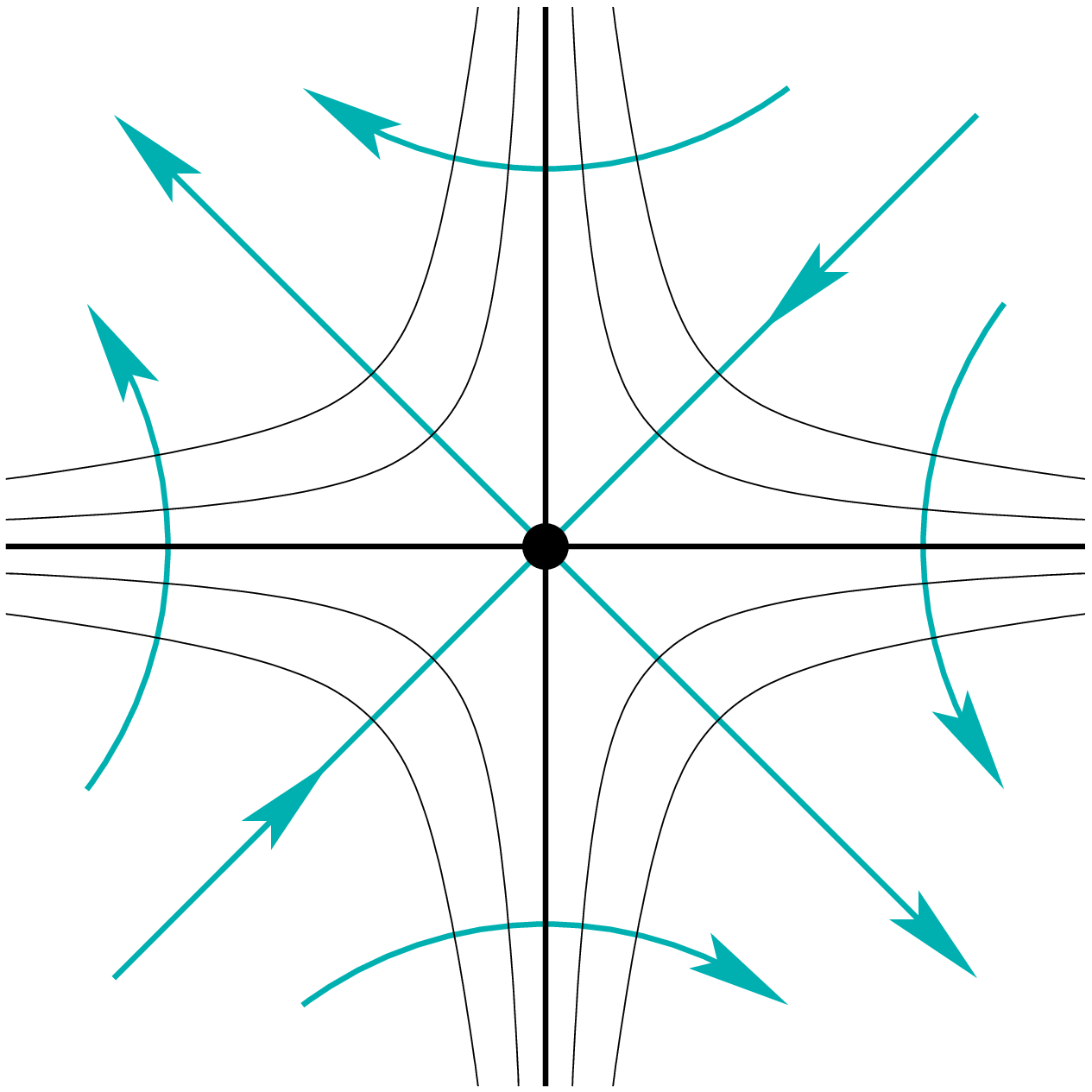}
 \refstepcounter{figure}\label{tranhyp}
 \\
 \centerline{\sc figure~\thefigure.}}
 \hskip 1.cm
 \parbox{5cm} {\includegraphics[scale=.46]{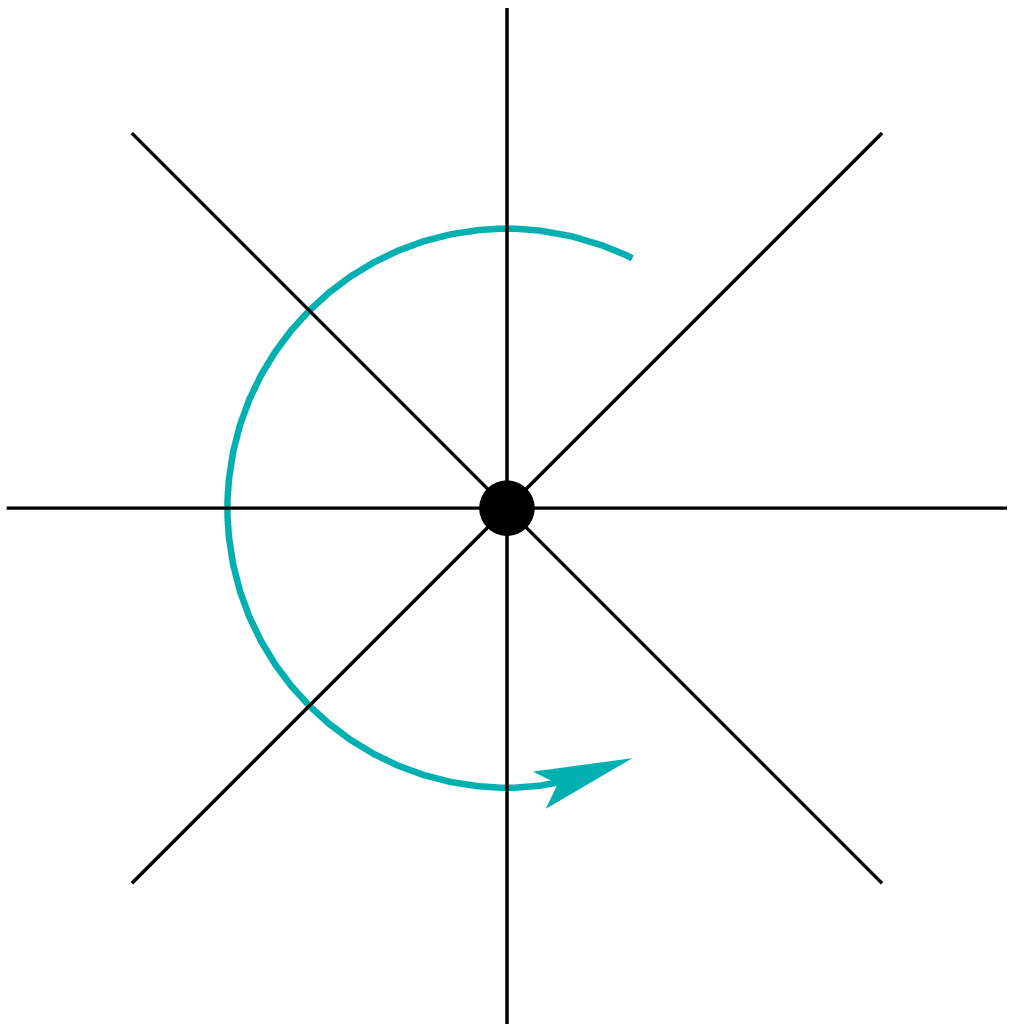}
 \refstepcounter{figure}\label{tranell}
 \\
 \centerline{\sc figure~\thefigure.}}
 }
 \vskip .5cm
 \centerline{
 \parbox{11cm}{\openup -2pt \small Figures~\ref{tranhyp} 
 and~\ref{tranell} show the behaviour
 of the foliation near a binding orbit $P$. The periodic orbit is shown
 perpendicular to the page and we only draw the trace of the foliation 
 on the page as black lines. 
 The big dot is the trace of the periodic orbit.
 \\
 \phantom{quad} In figure~\ref{tranhyp}
 the binding orbit has Conley-Zehnder index 2, and then it is hyperbolic 
 with positive eigenvalues. The inward pointing arrows 
 are the trace of the local stable manifold and the outward point arrows
 are the trace of the local unstable manifold.
 The black lines are the traces of 1-parameter families of leaves which 
 split into two rigid surfaces when they approach the periodic orbit.
 Two rigid surfaces are horizontal and the other two are vertical.
 The grey curves indicate the flow.
 \\
 \phantom{quad} In figure~\ref{tranell} the binding orbit 
 has Conley-Zehnder index 1 or 3.
 It is either elliptic or hyperbolic with negative eigenvalues.
 The foliation looks the same in both cases. When the orbit is hyperbolic
 the stable and unstable manifolds are M$\ddot{\text{\rm o}}$bius bands which 
 intersect transversally the leaves of the foliation.
 The grey curve represents the flow.
 }
 }
 
 \bigskip

 Consider the Reeb flow in a  neighbourhood of a binding periodic 
 orbit $P=(x,T)$ of index $\mu(P)\in\{1,3\}$. The $P$ is either elliptic 
 or $(-)$-hyperbolic. Fix a disk map $v_D:D\to S^3$ with 
 $v_D(e^{2\pi t})=x(tT)$ and a nowhere vanishing section $Z$ of
 the pull back $v_D^*\xi$ of the contact bundle $\xi$.
 The Reeb flow turns around the periodic 
 orbit at a speed  approaching the argument of the eigenvalues 
 of the linearized Poincar\'e map. When $\mu(P)=1$ and $P$ is 
 elliptic, the angular speed may be slow but bounded away from zero.
 The foliation gives an open book decomposition of a tubular 
 neighbourhood of $P$ and the flow is transversal to the foliation.
 All the transverse asymptotic vectors $e(t)$ of the leaves with asymptotic 
 limit $P$ must have the same the same (finite) winding number 
 with respect to $Z$. Then the orbits of the flow in a neighbourhood 
 of $P$ return (in the future and in the past) to each leaf 
 with asymptotic limit $P$ at a uniformly finite time
 (c.f.~\cite[lemma 5.2]{HWZ98}). 
 Since there are only a finite number of rigid surfaces, there
 must be at least one 1-dimensional family of leaves $F_\tau$
 with index $\Ind(F)=2$.
 If $f\la_0$ is not dynamically convex, 
 then the family decomposes on two rigid surfaces,
 at least one of them having $P$ as an asymptotic limit.
 Thus if  $f\la_0$ is not dynamically convex,  
 then on a neighbourhood of $P$, the orbits have uniformly finite 
 return times to a rigid surface.
 
 In the case of a non dynamically convex contact form,
 the points which lie in a stable manifold of a binding orbit $P$
 of $\mu$-index equal to 2 eventually do not return in the future to a 
 rigid surface. Similarly, the points in the unstable manifold of $P$ 
 eventually do not return in the past to a rigid surface.
 Any other point $x$ which is not in a binding periodic orbit
  returns to a rigid surface in a finite 
 (but not uniformly bounded) time. Because 
 flowing it a bit if necessary, it lies in an $I$-family 
 $F_\tau$ of leaves of index $\Ind(F_\tau)=2$.
 The parameter $\tau(x,t)$ defined by $\vr_t(x)\in F_{\tau(x,t)}$
 is strictly monotonous. Its derivative does not approach to zero
 if it is not in the local stable manifold of a binding orbit
 of $\mu$-index equal to 2. Then it reaches one end of the
 interval where the leaves decompose into rigid leaves.

 Hence in the non dynamically convex case:
 \begin{itemize}
 \item The set of rigid surfaces $\sec:=\cS_1$ 
       is composed of finitely many
       connected components, each with finite $d\la$-area and the form
       $d\la$ is non-degenerate on each of them.
 \item The points which are not in the stable manifold of a 
       binding periodic orbit of Conley-Zehnder index 2 return to
       $\sec$ infinitely many times in the future.
 \item The first return map $f$ is an area preserving
       local diffeomorphism. It is defined in all $\sec$ but
       the first intersection of the unstable manifolds of 
       binding periodic orbits of index 2 with $\sec$.
 \end{itemize}

  In the dynamically convex case,
   the foliation contains an $S^1$-family of leaves with index $\Ind(F)=2$.
  Then $S^3$ decomposes as the leaves of the
  family plus the binding orbits: $S^3=\cS_0\cup S_2$. Since the
  foliation contains a finite energy plane (with only one puncture),
  the leaves of the family are topological open discs and there is only
  one binding periodic orbit $P$ with index $\mu(P)=3$ by table~\eqref{tabla2}.
 This gives an open book decomposition of $S^3$.  In such case
  it has no rigid leaves but any leaf of the family is a (disc-like)
  global transversal section for the Reeb flow restricted to 
  $S^3\setminus\{P\}$. Thus in the dynamically convex case we have that
 \begin{itemize}
 \item every orbit of the flow
  intersects infinitely many times the (disc-like) section $\sec$, 
 \item $d\la\vert_\sec$ is an area form on $\sec$,
 \item the first return time is finite and uniformly bounded,
 \item the first return map $f:\sec\to\sec$ is an area preserving 
        diffeomorphism, with area form $d\la\vert_\sec$.
 \end{itemize}

 Define the diameter of a set $A\subset S^3$ by
 $$
 \diam(A):=\sup_{x,y\in A}d_{S^3}(x,y)
 $$

 We shall need the following

 \begin{proposition}\label{diamW}\quad

 Suppose that the Reeb vector field of the contact form $\la=f\la_0$
 is Kupka-Smale, i.e.  all its periodic orbits non-degenerate and
 the stable and unstable manifolds of its hyperbolic orbits, 
 when they intersect do it transversally.
 Furthermore, suppose that  $(S^3,\la)$
 is not dynamically convex.
  
  Then 
 there exists $A>0$ such that if $P$ is a binding orbit
 of Conley-Zehnder index $\mu(P)=2$ and $\cW$ is a connected component
 of $W^s(P)\cap\bS$ or $W^u(P)\cap\sec$, then
 $\diam \cW\ge A$.
 \end{proposition} 
 
 \begin{proof}\quad
 
 Let $\P$ be a binding periodic orbit of Conley-Zehnder index $\mu(\P)=2$. 
 We will prove the proposition for $W^s(\P)\cap\sec$.
 The proof for the unstable manifold is similar, flowing 
 the Reeb flow forwards.

 Since the Conley-Zehnder index of $\P$ is 2,  $\P$ is
 (+)-hyperbolic, i.e. its eigenvalues are positive.
 Then the stable manifold $W^s(\P)$ has two connected components
 which are immersed open cylinders.  
 Observe that $d\la\vert_{W^s(\P)}\equiv 0$, because
 $W^s(\P)$ is invariant under the Reeb flow and hence its 
 tangent space contains the Reeb vector field $X$, but
 $i_Xd\la\equiv 0$ on $S^3$.

 \centerline {\includegraphics[scale=.46]{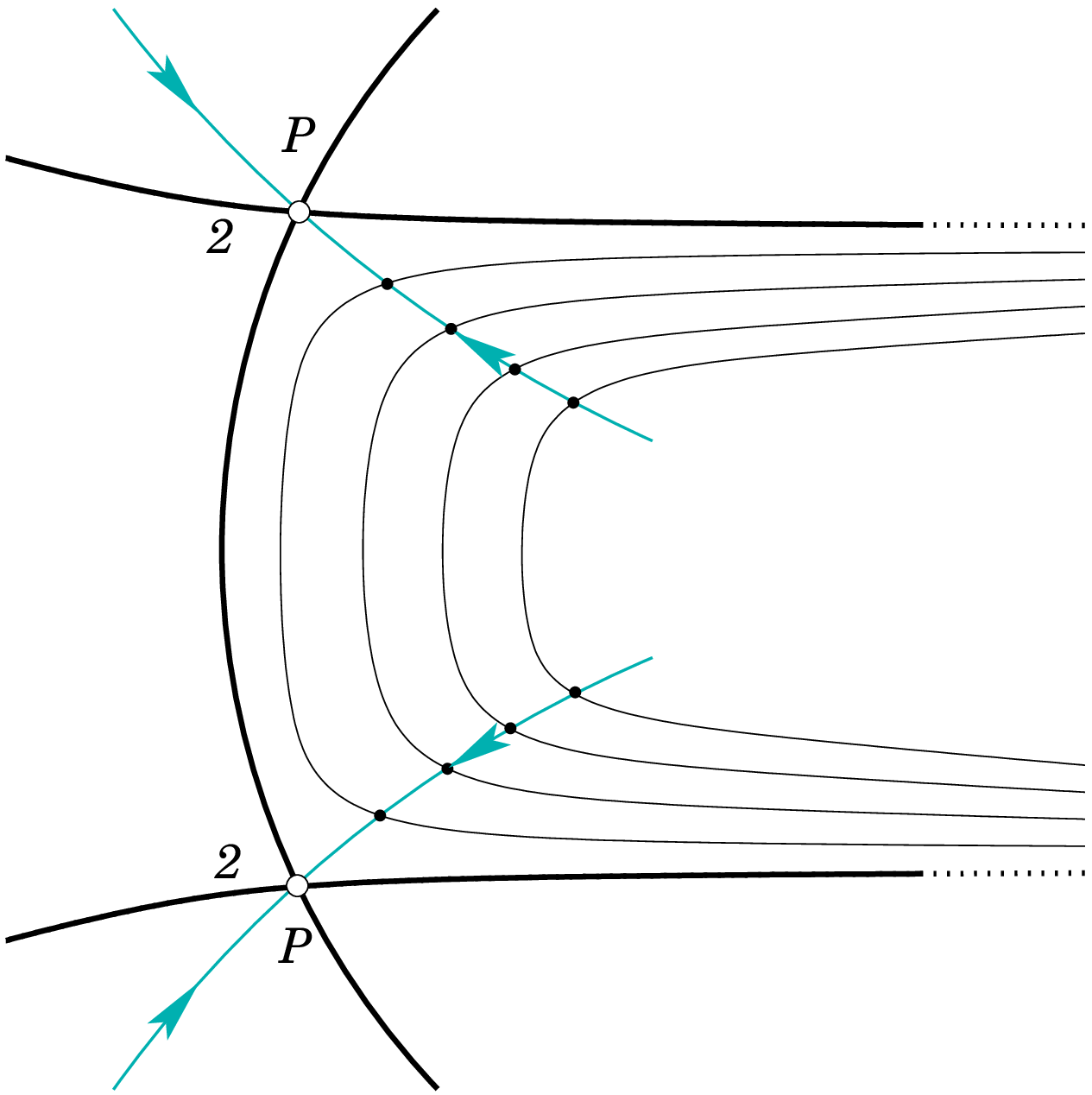}}
 \refstepcounter{figure}\label{loop}
 \vskip .5cm
 \centerline{
 \parbox{11cm}{\openup -2pt\small{\sc Figure~\thefigure.} This figure shows
 how the stable manifold of a binding orbit of Conley-Zehnder index 2
 cuts embedded circles on nearby leaves of the foliation belonging to
 a 1-parameter family.}}
 
 \medskip

 Choose one connected component of $W^s(\P)$.
 By lemma~7.6 in~\cite{HWZ1} there is a 1-dimensional 
 family of leaves $F_\tau$ parametrized by $\tau\in[0,1]$
 which at $\tau\to 0$ decomposes into two rigid surfaces $(C^+,C^-)$
 having $P$ as a common asymptotic limit such that the chosen component
 of the local stable manifold $W^s_{loc}(\P)$ intersects
 $F_\tau$ in a smooth embedded circle $S_\tau$ for all $\tau$ close to $0$.
 See figure~\ref{loop}. Flowing $W^s_{loc}(\P)$ backwards,
 since the flow is transversal to the foliation,
 the stable manifold keeps intersecting the leaves $F_\tau$
 in smooth embedded circles for all $\tau\in]0,1[$.
 When $\tau\to 1$, the family decomposes into two rigid orbits
 $(C_1^+,C_1^-)$ with a common binding orbit $P_1$ of index 2.  
 Recall that (lemma 5.2 in~\cite{HWZ98}) 
 in a neighbourhood of a binding orbit of odd $\mu$-index
 the solutions make a full turn around the binding orbit
 in a short interval of time. When the family of loops
 $\bS_\tau$ arrives at $\tau\to 1$ to the rigid surfaces
 either the connected component of intersection
 $W^s(\P)\cap(C_1^+\cup C_1^-)$ is inside one of the
 rigid surfaces, and in that case the connected component
 $\bS_1$ is also an embedded circle, or it intersects both surfaces.
 
 \begin{lemma}\label{diacir}
 There exist a positive constant 
 $$
 0<A<\min\{\,\diam(P)\;|\;P\text{\rm \; binding orbit}\,\},
 $$ 
 such that if
 \begin{itemize}
 \item[-] $R$ is a rigid surface of the foliation given by theorem~\ref{HWZ}.
 \item[-] $P$ is a binding orbit of Conley-Zehnder index $\mu(P)=2$.
 \item[-] $S$ is a connected component of $W^s(P)\cap R$ 
          or of $W^u(P)\cap R$ which is an embedded circle in $R$.
 \end{itemize}
 
 Then \quad $\diam(S)>A$.
 \end{lemma}

 \begin{proof}
 We only prove the lemma for the stable manifold, the same
 prove applies to the unstable manifold.
 Observe that the binding periodic orbit $P$ must be (+)-hyperbolic.
 Then the connected components of $W^s(P)$ are two topological cylinders.
 Let $\cW$ be the connected component of $W^s(P)$ which contains $S$.
 Then $S$ is also an embedded circle in $\cW$. Hence its homotopy 
 class is either 0 or $\pm 1$ in $\pi_1(\cW)=\Z$. 

 If $S$ is homotopically trivial then it bounds a disk in $\cW$. 
 Since the Reeb vector field $X$ is tangent to $\cW$ and has no singularities, 
  there is a point $s\in S$ where $X(s)$ is tangent to $S$.
 Since $S$ is also included in the leaf $R$ and the vector field $X$ is 
 transversal to $F$, $X(s)$ can not be tangent to $S$.
 Therefore $S$ is not homotopically trivial.

 Then $S$ is either homologous to $P$ or to $-P$ inside $\cW$.
 Choose the orientation in $S$ such that it is homologous to $P$.
 Since $S$ is a simple closed curve in $\cW$,  there is 
 an embedded surface $Q\subset \cW$ such that as a 2-chain 
 $\partial Q= S-P$. Since $d\la\vert_{W^s(P)}\equiv 0$, by
 Stokes theorem
 $$
 \oint_S\la=\int_Q d\la+\oint_P\la=\oint_P\la=T(P),
 $$
 where $T(P)$ is the minimal period of $P$ (recall that the
 binding orbits are simply covered).

 Let $\Y$ be the set of closed continuous curves $y:S^1\to\sec$ 
 inside a rigid surface which are not homotopically trivial.
 Since the set of rigid surfaces is finite and their ends 
 are a finite set of closed orbits, 
 $$
 a:=\tfrac 12\inf\{\,\diam(y)\,|\,y\in\Y\,\}>0
 $$
 is positive. Then if $S$ is not homotopically trivial on
 the rigid surface $R$ then
 $$
 \diam(S)>a.
 $$ 

 If $S$ is homotopically trivial on $R$ then it bounds an embedded disk 
 $D\subset R$. Then
 $$
 \int_D d\la=\oint_S\la>T(P).
 $$Let 
 $$
 b:=\min\{\,T(P)\,\vert\, P\text{ binding orbit }\,\}>0.
 $$
 Since the rigid surfaces are embedded, transversal to the 
 non-singular vector field $X$, with compact closure
 and finitely many, 
 $$
 c:=\inf\{\,\diam D\,\vert\,D\hookrightarrow\sec \text{ embedded disk},
 \;d\la\text{-area}(D)>b\,\}>0
 $$
 is positive. Then 
 $$
 \diam(S)=\diam(D)>c.
 $$
 Now let
 $$
 A:=\min\{\,a,\;b,\;c\,\}>0.
 $$
 \end{proof}

 In the case the family $\bS_\tau$ intersects at an embedded circle
 $\bS_1$ only one of the rigid surfaces, say $C_1$, we use 
 lemma~\ref{diacir} to show that $\diam(\bS_1)>A$ and continue flowing 
 backwards along a family of leaves $F_\tau$, $\tau\in]1,2[$.
 The intersection $\bS_1$ continues as  embedded circles $\bS_\tau$
 for all $\tau\in]1.2[$ and we repeat the argument.

 \parbox{7cm}{
 \centerline{\includegraphics[scale=.46]{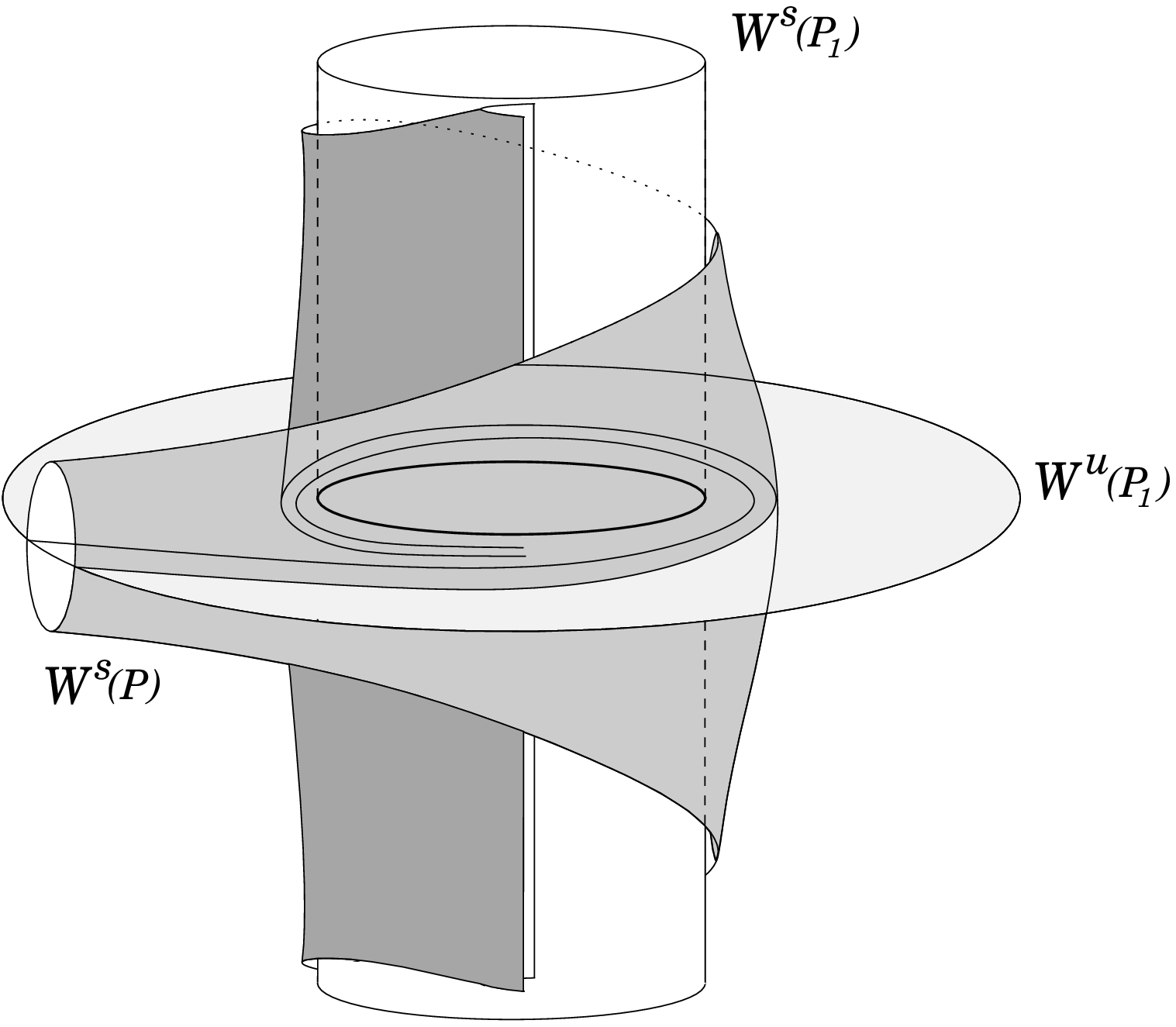}}
 \vskip 10pt
 \refstepcounter{figure}
 %\centerline{{\sc figure \thefigure\label{intf}\phantom{aaa\,}}}
 \vskip .8cm
 \openup -2pt
 \vskip -22pt
 {\sc figure \thefigure.\label{int}}
 { \small This figure shows how the
 stable manifold of an hyperbolic periodic orbit $P$
 accumulates on the whole stable manifold of another
 hyperbolic periodic orbit when there is an heteroclinic
 transversal intersection.
 \phantom{quiero una linea en blanco tengo que po  } }
 }
 \hskip .5cm
 \parbox{7cm}{
 \centerline {\includegraphics[scale=.46]{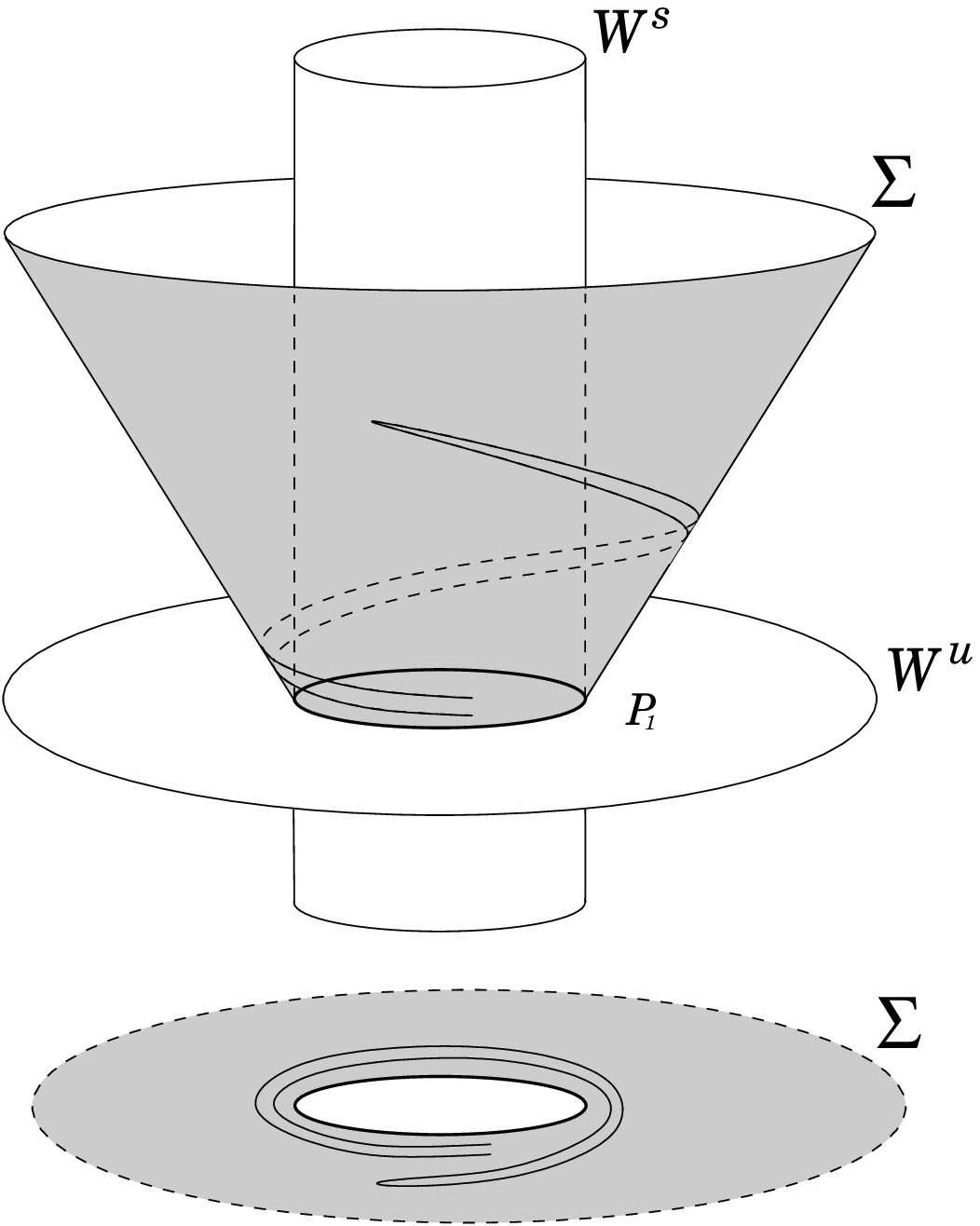}}
 \refstepcounter{figure}
 %\centerline{\sc figure \thefigure\label{sec}\phantom{a\;}}
 \vskip .8cm
 \openup -2pt
 {\sc figure~\thefigure\label{sec}}.
 {\small This figure shows the approach of  rigid surface $\Si$ 
  to an asymptotic limit $P_1$ of index $\mu(P_1)=2$.
  If the unstable manifold of $P_1$ intersects transversally
  the stable manifold $W^s(P)$ of another binding orbit $P$,
   as in figure~\ref{int},
  then figure~\ref{sec} also shows
  how the \; {in-} 
   }}
  \vskip -10pt
  \phantom{aa}\parbox{14.6cm}{\openup -2pt
   {\small \,\,
   tersection of $W^s(P)$
  with the rigid surface $\Si$ accumulates on the whole asymptotic
  limit $P_1$.
   The lower part of the figure is a representation of a
   neighbourhood
  of the puncture corresponding to $P_1$ in the rigid surface
  $\Si$ and the intersection $W^s(P)\cap \Si$.
  }}

  \medskip

 In the case the family $\bS_\tau$ hits both rigid surfaces 
 $C_1^+$ and $C_1^-$ then in the limit $\tau\to 1$ the stable 
 manifold must intersect the unstable manifold of 
 binding orbit $P_1$ of index 2 which is a common asymptotic limit
 to $C_1^+$ and $C_1^-$. The family $S_\tau$ decomposes into three 
 sets: a subset in $W^s(\P)\cap W^u(P_1)$, which for negative time
 has no return to the union of rigid surfaces $\sec$;
 and two connected components of $W^s(\P)\cap\sec$, one component
 $\bS_1$ in $C_1^+$ and another $\bS_1^-$in $C_1^-$. 
 Both components are immersions 
 of the real line $\re$ in $C_1$ or $C_1^-$ and both contain
 the periodic orbit $P_1$ in their closure. Hence
 $$
 \diam(\bS_1^+)\ge\diam(P_1)>A,
 $$
 and similarly $\diam(\bS_1^-)>A$.
 
 Choose one of the components, say $\bS_1^+$. When we continue flowing 
 backwards $\bS_1^+$, it accumulates on a branch of the stable manifold 
 $W^s(P_1)$ by the $\la$-lemma. The argument repeats, the branch of 
 $W^s(P_1)$ cuts a family of surfaces $F_\tau$ parametrized by 
 $\tau\in]1,2[$ on embedded circles $S^1_\tau$, and the backward flow 
 of the component $\bS^+_1$ intersects the family on connected 
 subsets which accumulate on {\it the whole} 
 $S^1_\tau\subset W^s(P_1)\cap F_\tau$, (see
 figures~\ref{int},~\ref{sec}).

 When $\tau\to 2$ the family $F_\tau$ decomposes into two rigid surfaces 
 $C_2^+$, $C_1^-$ and the family of connected components 
 $\bS^+_\tau\subset W^s(\P)\cap F_\tau$ either intersects only 
 one of them or both. If it intersects only one, say $C_2^+$, 
 the intersection $\bS_2$ is in the interior of $C_2^+$ 
 and its closure contains the intersection
 $S_{\tau=2}^1\subset W^s(P_1)\cap C_2^+$, which is an embedded circle.
 Then
 $$
 \diam(\bS_{\tau=2}^+)\ge \diam(S_{\tau=2}^1)>A.
 $$

 \centerline {\includegraphics[scale=.46]{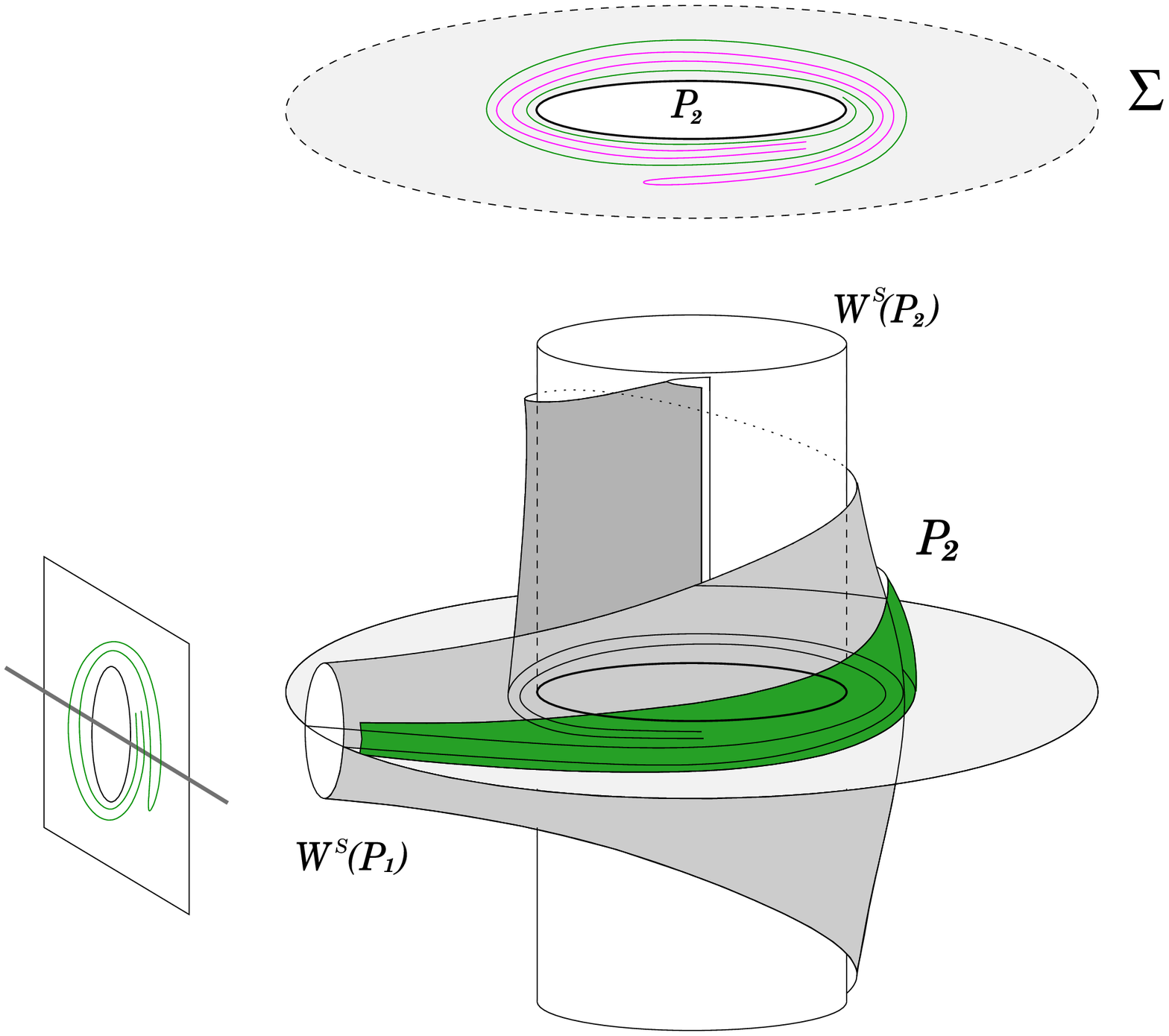}}
 \refstepcounter{figure}\label{recu}
 \vskip .2cm
 \centerline{\parbox{11cm}{\openup -2pt \small{\sc figure~\thefigure.}
 This figure shows how the backward flow of a small interval $J$ 
 of a stable manifold $W^s(P)$
 which intersects transversally the unstable manifold of a binding orbit
 $P_2$ accumulates on the unstable manifold of $P_2$ and its intersection
 on a rigid surface with asymptotic limit $P_2$ contains the orbit $P_2$
 in its closure. In this case, $W^s(P)$ first accumulates on $W^s(P_1)$
 which, in turn, accumulates on $W^s(P_2)$.}}
 
 \bigskip

 If the family $\bS^+_\tau$ as $\tau\to 2$ intersects both surfaces
 $C^+_2$ and $C^+_1$, since $\bS_\tau^+$ is connected, it intersects
 (transversally) $W^u(P_2)$. Choose for example $C^+_2$. 
 Let $\bG$ be a connected component of the intersection of 
 $\bS_\tau^+$ as $\tau\to 2$ with $C_2^+$. The $\bG$ is a 
 1-dimensional submanifold whose forward flow contains an open segment
 $\bJ\subset W^s(\P)$ with an endpoint in $W^u(P_2)$. 
 which is transversal to $W^u(P_2)$. Flowing backwards $\bG$ 
 as in figure~\ref{recu}, we have that the connected component 
 $\bG\subset C_2^+$ must accumulate on the whole periodic orbit $P_2$.
 Hence
 $$
\diam(\bG)\ge \diam(P_2)>A.
 $$
 Flowing backwards, we have that  $\bG$ accumulates on
 the unstable manifold $W^u(P_2)$, and the argument goes on.

 Alternatively, we could say that the connected component $\bG$
 accumulates on a whole component $B$ of $W^s(P_1)\cap C_2^+$, 
 which in turn accumulates on the whole periodic orbit $P_2$.
 Thus
 $$
 \diam(\bG)\ge\diam(B)\ge\diam(P_2)>A;
 $$
 and the argument goes on.
 \end{proof}

 \bigskip

  We shall be interested in returns to a rigid surface.
  Let $\Si$ be a rigid surface. Define the return times
  $\tau_n:\Si\to\re\cup\{+\infty\}$ and return maps
  $F_n:[\tau_n<+\infty]\to\Si$ to $\Si$ by
  \begin{align*}
  \tau_0 &:\equiv 0, 
   \\
   F_n(x) 
    &:= \vr_{\tau_n(x)}(x), \qquad \text{ if }\quad\tau_n(x)<+\infty\,; 
    \\
   \tau_{n+1}(x) 
    &:=\inf\{\,t>\tau_n(x)\;|\; \vr_t(x)\in\Si\,\}.    
  \end{align*}

   \bigskip

  \begin{proposition}\label{diamF}
  Suppose that the pair $(S^3,\la)$ satisfies the hypothesis 
  of proposition~\ref{diamW}.
  Let $A>0$ be the constant given by proposition~\ref{diamW}.
  Let $B$ be an open subset of a rigid surface $\Si$
  such that $\tau_N\vert_B$ is finite but  not bounded
  for some $N>0$.
  Then $\diam(F_N(B))>A$.
  \end{proposition}
  
  \begin{proof}
  Since $\tau_N\vert_B$ is not bounded, there is a point $q$ in the
  boundary $\partial B$ such that $\tau_N(q)=+\infty$. Let
  $M:=\max\{\, k\ge 0\;|\; \tau_k(q)<+\infty\,\}$.
  Since $\tau_0\equiv 0$, we have that $0\le M\le N-1$. 

 %  Since $\tau_N(I)\equiv+\infty$, then the points $x\in I$ accumulate on
 %  $P$ before they have an $N$-th return to $\Si$. For $x\in I$ define
 %  $$
 %  M(x):=\max\,\{\,n\ge 0\;|\;\tau_n(x)<+\infty\,\}.
 %  $$
 %  For $0\le k<N$, the sets $[M=k]\cap I$ are either empty or
 %  subintervals of $I$. Since $I$ is a finite union
 %  $I=\cup_{k=1}^{N-1}([M=k]\cap I)$, reducing $I$ if necessary, we can
 %  assume that $I\subseteq [M=k_0]$ and we shall write
 %  $M(x)=M:=k_0$.

 \vskip .8cm

 \includegraphics[scale=.46]{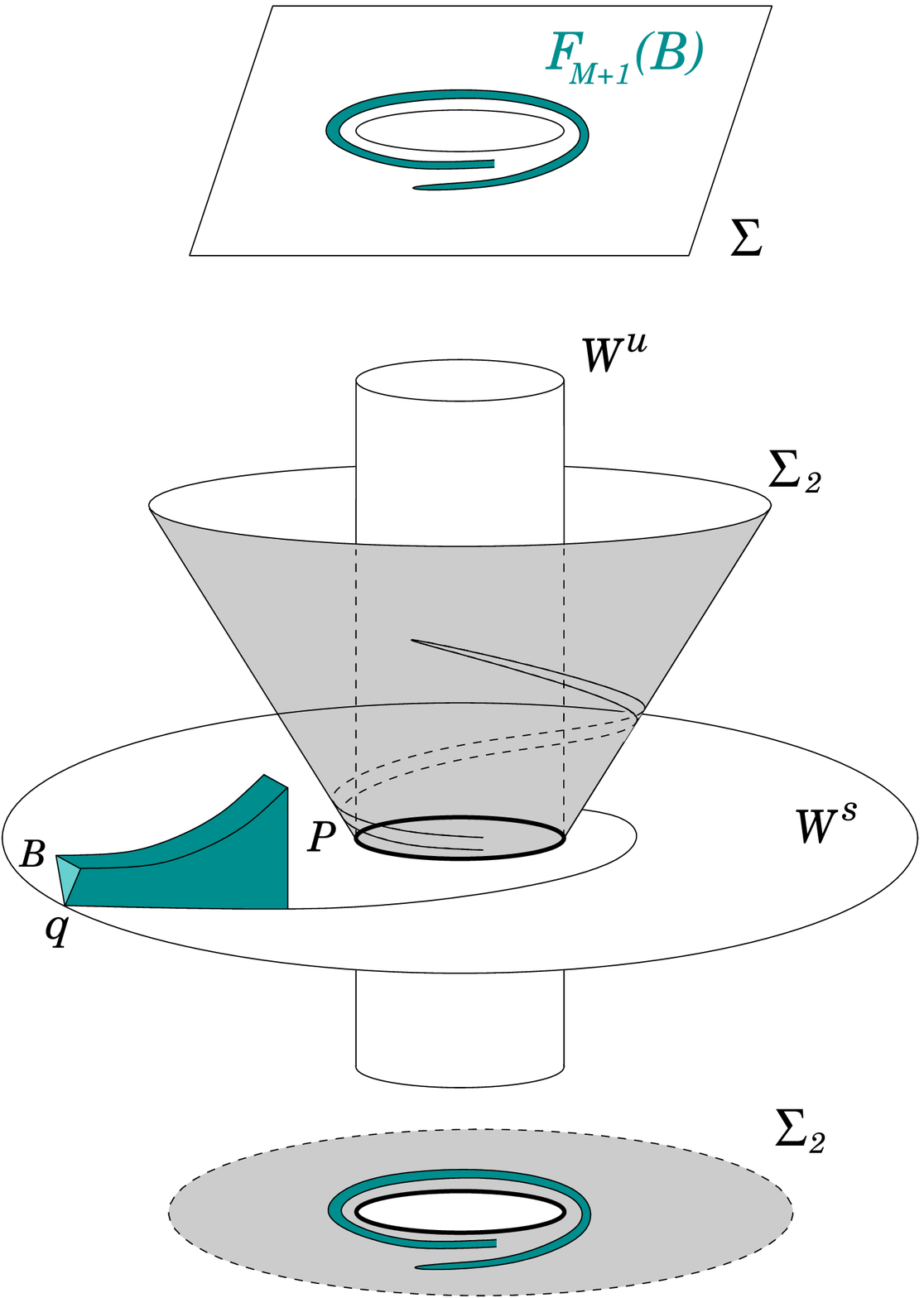}
 \refstepcounter{figure}\label{B}
 \qquad\qquad
 \vskip -8cm \hskip 8cm \includegraphics[scale=.46]{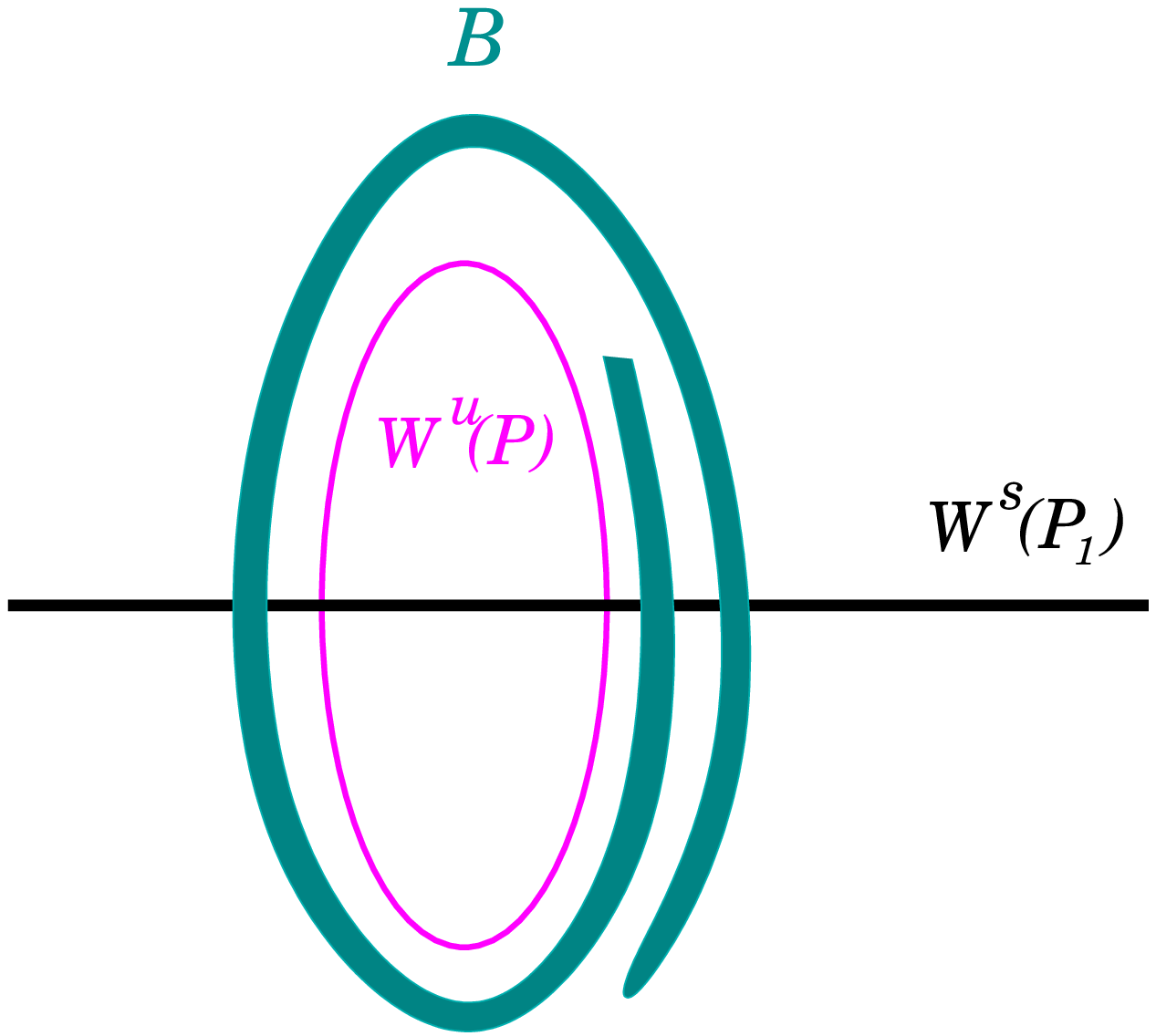}
 \vskip .5cm \hskip 7cm
 \refstepcounter{figure}\label{Bint}
 \parbox{7cm}{\openup -2pt \small{\sc figure~\thefigure.}
 This figure shows how if $W^u(P)$ intersects another stable
 manifold $W^s(P_1)$ before returning to $\Si$, then
 $W^s(P_1)$ must intersect the interior of the
 forward iterate of $B$ before it returns to $\Si$. }

 \vskip .4cm \hskip 1cm
 \parbox{13cm}{\openup -2pt \small{\sc figure~\ref{B}.}
 The figure shows how the forward flow of  an open set $B$
 in a rigid surface $\Si$ 
 with a boundary point $q$ in the  stable manifold $W^s(P)$
 intersects a rigid surface with asymptotic limit $P$
 in a subset which accumulates on the {\it whole} orbit $P$.
 If its forward flow does not intersect other stable
 manifolds of binding perioidc orbits of index 2, then
 its returns to $\Si$ are sets which contain in their closure
  an intersection of the unstable manifold $W^u(P)$ 
  of $P$ which is a topological circle.}
 \vskip .4cm

  Since $M<N$,  $\tau_M|B<+\infty$. By the continuity of the flow
  the set $[\tau_M<+\infty]$ is open.  
   The map $F_M$ is a local
  diffeomorphism on the open set $[\tau_M<+\infty]$. 
  We have that $F_M(B)\subseteq\Si$,
  $F_M(q)\in\partial F_M(B)$. Moreover,
  $F_M(q)\in W^s(P)\cap \Si$ for a binding
  periodic orbit $P$ with index $ \mu(P)=2$.
  Flowing forward $F_M(B)$  it reaches a
  neighbourhood of $P$, then intersects a rigid surface accumulating
  on the whole orbit $P$ as in figure~\ref{B}, and then it follows
  accumulating on $W^u(P)$. The forward orbit of $F_M(q)$ accumulates
  on the periodic orbit $P$ and does not return to $\Si$. 

  Suppose that  $W^u(P)$ intersects another
  stable manifold $W^s(P_1)$ of a binding periodic orbit $P_1$ with
  $\mu(P_1)=2$, with point accumulating on $P_1$ 
  before it intersects $\Si$.
  Then by hypothesis $W^u(P)$ intersects $W^s(P_1)$
  transversally. Since the forward flowing of $F_M(B)$
  accumulates on the whole $W^u(P)$, it intersects
  $W^s(P_1)$ before it returns to $\Si$ (see figure~\ref{Bint}).
  This implies that there are points  $x\in B$ for which
  $\tau_{M+1}(x)=+\infty$. Since $M+1\le N$, this contradicts 
  the hypothesis $\tau_N|_B<+\infty$.

  Then the unstable manifold $W^u(P)$ intersects 
  $\Si$ in a circle $S_{M+1}$ and $F_{M+1}(B)$ accumulates
  on the whole circle $S_{M+1}$. By lemma~\ref{diacir},
  $$
  \diam(F_{M+1}(B))=\diam(\ov{F_{M+1}(B)})\ge \diam S_{M+1}>A.
  $$
  If $M+1<N$, for each of the following iterates $F_k$, $M+1<k\le N$,
  the argument above show that the forward flowing of $S_{M+1}$ 
  does not intersect a stable manifold of a binding orbit of index 2
  before returning to $\Si$. Then $S_{M+1}\subset [\tau_{N-M-1}<+\infty]$,
  $S_N:=F_{N-M-1}(S_{M-1})$ is a circle and $F_B(B)$ accumulates on
  $S_N$. Then $\diam F_N(B)\ge \diam S_N>A$.

  \end{proof}

 \medskip

\section{Proof of theorem~\refteo.}\label{proofs}

 The Kupka-Smale theorem for geodesic flows is proven in~\cite{CP2},
 we recall it here.
 Let $J_{s}^{r}(n-1)$ be the set of $r$-jets of symplectic 
 automorphisms of $\re^{n-1}\oplus \re^{n-1}$ that fix the origin.
 Clearly one can identify $J_{s}^{1}(n-1)$ with $Sp(n-1)$.
  A set $Q\subset J_{s}^{r}(n-1)$
 is said to be {\it invariant} if for all 
 $\sigma\in J_{s}^{r}(n-1)$, $\sigma Q\sigma^{-1}=Q$.
 
 \begin{Kupka}\quad\label{KS}
   Let $Q\subset J^{r-1}_s(n-1)$ be open, dense and invariant.
 Then there exists a residual subset $\cO\subset\cG^r$
 such that for all $g\in\cO$:
 \begin{itemize}
 \item  The $(r-1)$-jet of the Poincar\'e map of every
  closed geodesic of $g$  belongs to $Q$.
  \item All heteroclinic points of hyperbolic closed geodesics
  of $g$ are transversal.
 \end{itemize}
 \end{Kupka}

 % We  say that a set $\cH$ of closed geodesics for a riemannian manifold 
 % $(M,g_0)$ is {\it $C^2$ stably hyperbolic} if there exists a 
 % $C^2$ open neighborhood $\cU$ of $g_0$ in the set of $C^\infty$ 
 % riemannian metrics of $M$ such that for all $g\in\cU$, all the 
 % continuations of the periodic orbits in $\cH$ for the geodesic 
 % flow of $g$ are hyperbolic (i.e. their linearized Poincare map
 % has no eigenvalues of modulus 1).  

 \medskip

  Let $M$ be a closed 2-dimensional smooth manifold
  and let $\cH^1(M)$ be the set of $C^r$~riemannian metrics,
  $r\ge 4$ on $M$ all of whose closed geodesics are hyperbolic,
  endowed with the $C^2$ topology and let
  $\cF^1(M)=\inte\bigl(\cH^1(M)\bigr)$ be the interior of
  $\cH^1(M)$ in the $C^2$ topology.
  %Theorem~\ref{herman} says that $\cF^1(S^2)$ is empty.
  If $g\in\cF^1(M)$, let $\oPg$ be the closure of the periodic
  orbits of the geodesic flow of $g$.

    A hyperbolic set $\La$ is said {\it locally maximal}
  if there exists an open neighbourhood $U$ of $\La$, such that
  $\La$ is the maximal invariant subset of $U$, i.e.
  $
  \La = \bigcap_{t\in\re}d\phi_t^g(U).
  $
  A {\it basic set} is a locally maximal hyperbolic set with
  a dense orbit.
  It is {\it non-trivial} if it is not  a single
  closed orbit.
  It follows that in a basic set the periodic orbits are
  dense in its relative non-wandering
  set~\cite[cor. 6.4.19]{HK},~\cite{Bowen6,Bowen3}.

 \bigskip

\begin{theorem}[{\cite[prop. 5.5, cor. 5.9]{CP2}}]\label{gpe}
\quad

If $g\in\cF^1(M)$ then the closure $\oPg$ is hyperbolic.
Moreover, $\oPg$ decomposes into a finite number of hyperbolic basic
sets and at least one of them is non-trivial.

\end{theorem}

 The non-triviality of the basic set in this theorem
 (c.f. \cite[cor. 5.9]{CP2}) is obtained from the fact that
 the geodesic flow for $M$ has infinitely many closed orbits.
 If $M$ is non simply connected, then $\pi_1(M)$ is infinite, 
 and any riemannian metric has a minimizing closed geodesic in
 each free homotopy class. If $M=S^2$, then Bangert~\cite{Bang0}
  and Franks~\cite{Fra2} proved that $M$ has infinitely many
  geometrically distinct geodesics. For the projective space
  $\bP^2$, lifting the riemannian metric on $\bP^2$ to its double cover
  $S^2\to\bP^2$ one obtains that also $\bP^2$ has infinitely many
  closed geodesics.
 
  \medskip

 The arguments in the above theorem need $C^2$ perturbations.
 For example, it is not known if a riemannian 
 metric
 whose geodesic flow is topologically equivalent to an Anosov flow
 can be $C^r$ approximated, $r\ge 3$, by a metric whose 
 geodesic flow is Anosov.

  \medskip
  
   Define an equivalence  relation on $\Pg$ by saying that $x\sim y$
   if and only if $W^s(x)\cap W^y(y)\ne\nulo$ and $W^u(x)\cap
   W^s(y)\ne\nulo$. Its equivalence classes are called {\it homoclinic
   classes}. The spectral decomposition theorem~\cite{PS}, \cite[ex.
   18.3.7]{HK},
    used to obtain the basic set in theorem~\ref{gpe}, 
    states that the basic set is a {\sl whole homoclinic class}.

\bigskip

 \noindent{\bf Proof of theorem~\refteo:}
 
 Assume by contradiction that $(S^2,g_0)$ is a $C^\infty$  
 riemannian metric which can not 
 be $C^ 2$~approximated by one having an elliptic closed geodesic.
 Then $g\in\cF^1(S^2)$. We can assume that $g$ is a Kupka-Smale metric.
 By theorem~\ref{gpe}, there
 is a hyperbolic non-trivial 
 basic set $\La$ for the geodesic flow $\phi_t$ of $(S^2,g)$.
 Let $m$ be the normalized Liouville's measure
 for $(T^1S^2,\phi_t)$.
 
 Since $g$ is $C^r$, $r\ge 3$, we have that $f$ is $C^2$. 
 Bowen and Ruelle~\cite[cor. 5.7.(b)]{BoRu} 
 proved that since $f$ is $C^2$, if $\La$ 
 has positive measure, then $\La$ must be open. Since $\La$ is 
  closed, either $m(\La)=0$ or $\La=\Sigma$.
 If $\La=\Sigma$, then the geodesic flow must be Anosov.
 But there are no Anosov\footnote{because
 $\pi_1(T^1 S^2)$ would have exponential growth (Margulis~\cite{Margulis},
 Plante \& Thurston~\cite{PlanteThurston} or
 Klingenberg~\cite{Klingenberg}) or 
 also because $(S^2,g)$ would not have conjugate points 
 (Klingenberg~\cite{Klingenberg}).} geodesic flows for $S^2$.
 Thus $m(\La)=0$.
 
 Let 
 $
 Q=W^s(\La)\cup W^u(\La)
 =\{\,x\in\Sigma\,\vert\,\alpha\text{\rm -lim}(x)\subseteq \La
 \text{ or }\omega\text{\rm -lim}(x)\subseteq\La\,\}.
 $
 Then Hirsch, Palis, Pugh and Shub~\cite{HPPS}, proved that
 \begin{equation}\label{hpps}
 Q=\bigcup_{t\in\re} \bigcup_{x\in\La} \phi_t\big(W^{ss}_\e(x)\cup
 W^{uu}_\e(x)\big);
 \end{equation}
 where  for $x\in\La$,  $W^{ss}_\e(x):=\cap_{t\ge
 0}\,\phi_{-t}\big(B_\e(\phi_t(x))\big)$
 and $W^u_\e(x):=\cap_{t\ge 0}\,\phi_{t}\big(B_\e(\phi_{-t}(x))\big)$ 
 are the local strong stable and unstable
 manifolds.
 Since the homoclinic class $\La$ for the {\it flow} $\phi_t$ is transitive, $Q$ is
 connected\footnote{but its intersections with a transversal 
 ``surface of section'' may be not connected.}.
 The points in $Q\setminus\La$ 
 are wandering. By Poincar\'e's recurrence theorem  $m(Q\setminus\La)=0$.
 Hence $m(Q)=0$.

 We say that a point $x\in\La$ is an {\it interior point of $\La$}, if
 $x$ is an accumulation point of $\La$ on both connected components
 of $W^{ss}_\e(x)$ and also of $W^{uu}_\e(x)$. Since $\La$ is
 non-trivial, it has infinitely many interior points.
 
 We shall need the following two lemmas:
 
 \newpage
 
 \begin{Lemma}\label{homclass} \quad
 $\La=W^s(\La)\cap W^u(\La)$.
 \end{Lemma}
 
 \begin{proof}
 Clearly $\La\subseteq W^s(\La)\cap W^u(\La)$. Suppose that $z\in
 W^s(\La)\cap W^u(\La)$. By~\eqref{hpps} there are $x,y\in \La$ such
 that $z\in W^s(x)\cap W^u(y)$  or  $z\in W^u(x)\cap W^s(y)$. Suppose
 it is the first case, the second is similar.
 We have to prove that $z\in\La$.
 
 The point $z$ may be in a tangential intersection of $W^s(x)$ with
 $W^u(y)$ or in a topologically transversal intersection. The
 topologically transversal case is easier and uses the same methods, so
 we prove only the tangential case.

 \begin{figure}
 \includegraphics[scale=.5]{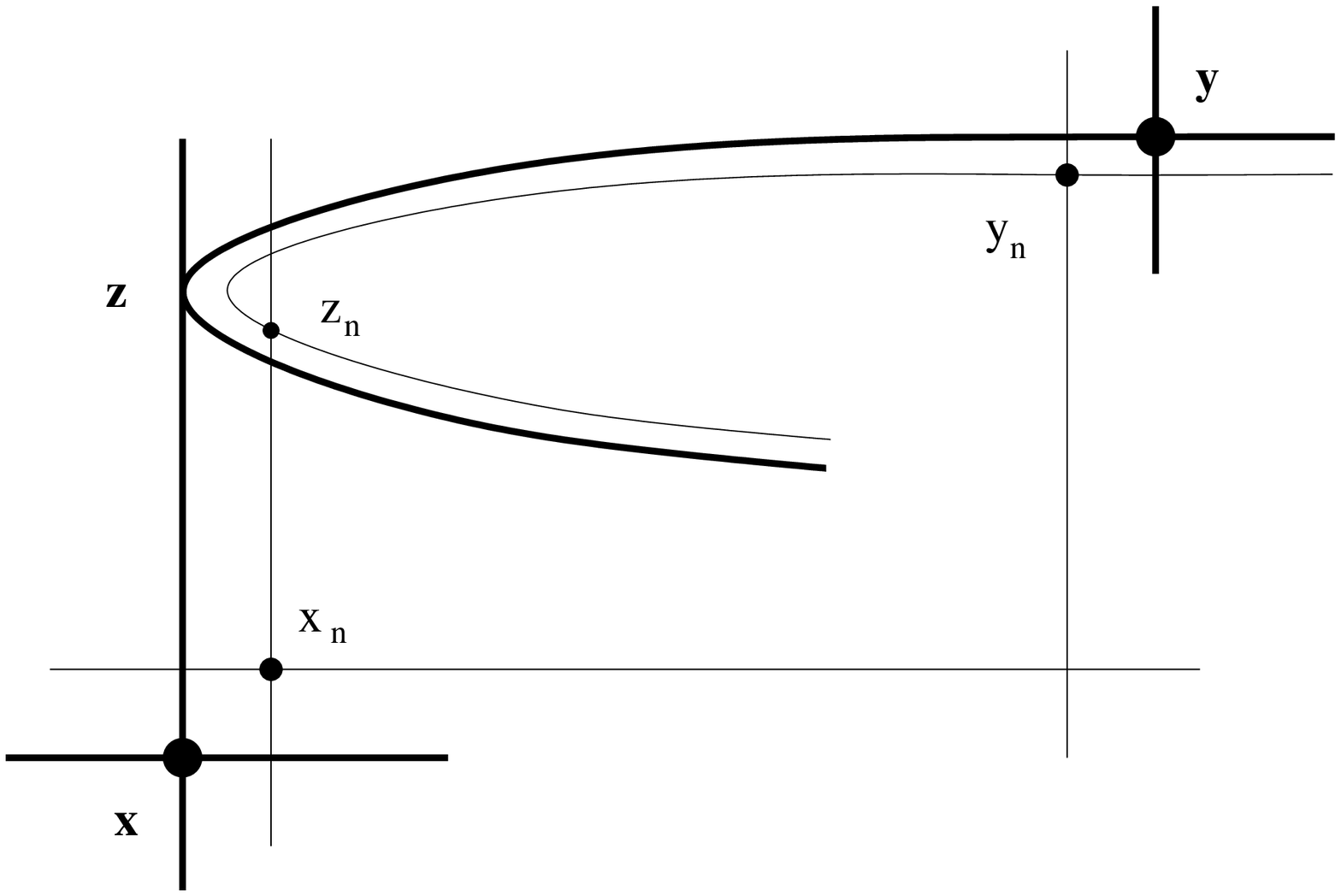}
 \caption{\small This figure shows how if one of $x,\,y\in \La$
 is an interior point of the homoclinic class $\La$, then any
 $z\in W^s(x)\cap W^u(y)$ is accumulated by homoclinic points $z_n$
 of periodic points $x_n, \; y_n$ in the same homoclinic class
 $\La$.}\label{tang}
 \end{figure}

 Suppose that $x$ is an interior point of $\La$. The case when
 $y$ is an interior point of $\La$ is similar. Then $x$ is accumulated 
 by points of $\La$ in its four quadrants made by $W^s_{loc}(x)$ and
 $W^u_{loc}(x)$. Since the periodic points are dense in $\La$, 
 there exist sequences $x_n\to x$, $y_n\to y$ in $\La\cap\Pg$
 and points $z_n\in W^s(x_n)\cap W^u(y_n)$ such that $\lim_n z_n=z$,
 as in figure~\ref{tang}.
 Since $g$ is Kupka-Smale and $x_n$, $y_n$ are periodic points,
 the intersections $W^s(x_n)\transv W^u(y_n)$ 
 are transversal. Since $x_n$, $y_n$ are in the homoclinic class $\La$,
  we also have that $W^u(x_n)\cap W^s(y_n)\ne\emptyset$. Then there
 are  hyperbolic sets $H$ inside the homoclinic class $\La$ which
 accumulate on $z_n$, see figure~\ref{horse}. Therefore $z_n\in\La$. 
 Since $\La$ is closed,  $z=\lim_n z_n\in\La$.
 
 \begin{figure}
 \includegraphics[scale=.5]{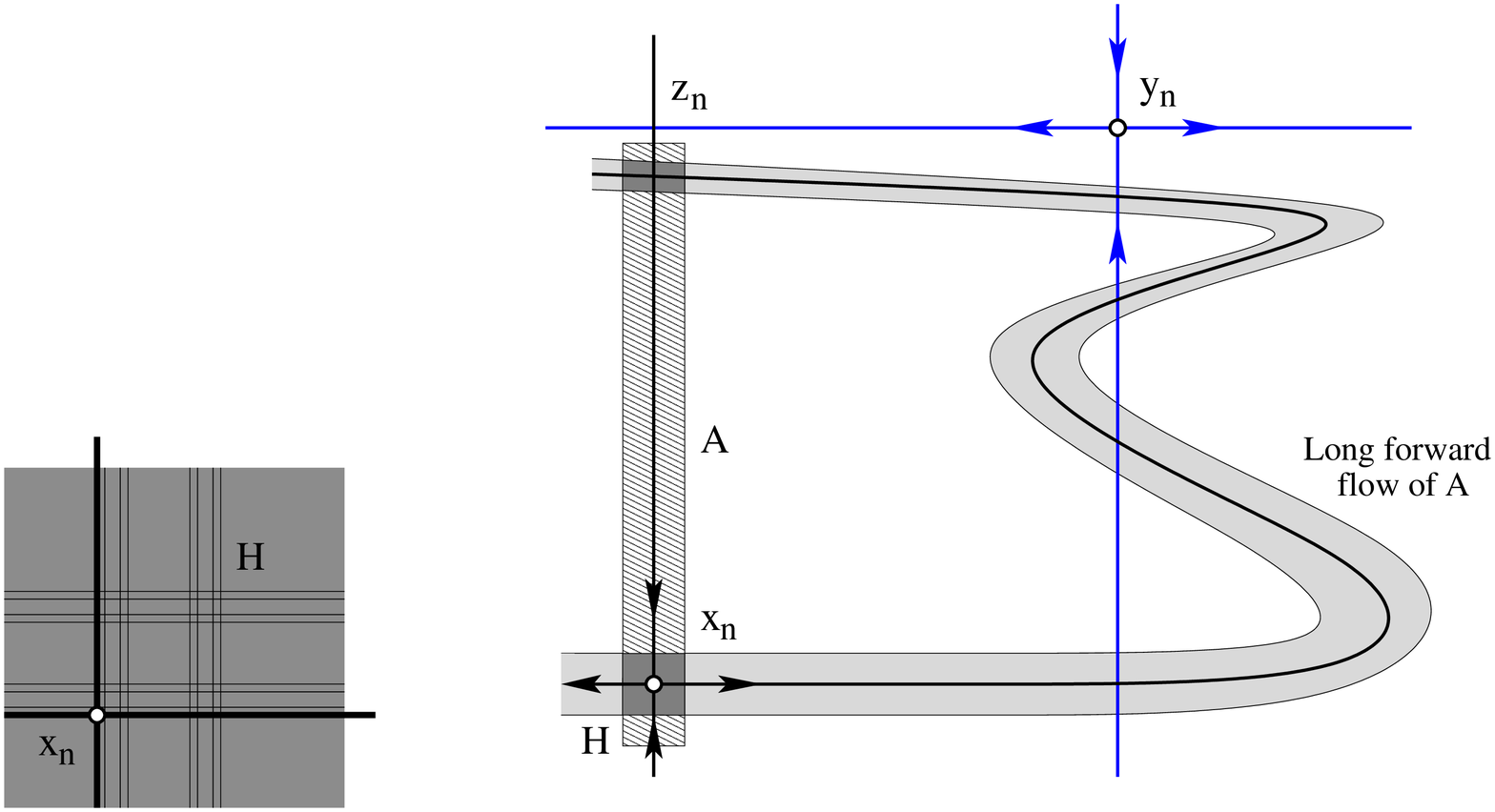}
 \caption{\small This figure shows why the point $z_n\in W^s(x_n)\cap
 W^u(y_n)$ is in the homoclinic class $\La$. Since also $W^u(x_n)\cap
 W^s(y_n)\ne\emptyset$, the point $z_n$ is accumulated by small 
 (transitive) horseshoes $H=H_m$, each of them also accumulating on $x_n$.
 By the local product structure on the uniformly hiperbolic set $\oPg$
 the homoclinic classes are (relatively) open subsets of $\oPg$. Thus, 
 when this construction is made with a sufficiently thin
 rectangle $A$ and a very large forward iterate by the flow, 
 the periodic points in the horseshoe $H$ 
 are in the same homoclinic class $\La$ as $x_n$.
 In fact the whole (transitive) horseshoe is inside the homoclinic
 class $\La$, because, by the
 shadowing lemma, for any $\e>0$, 
 there are $\e$-dense periodic orbits on $\La$. }\label{horse}
 \end{figure}

 Now suppose that $x$ and $y$ are not interior points of $\La$.
 Then (cf. ~\cite[prop. 2.1.1]{BoLa} or~\cite{NePa})  they
 are in the invariant manifolds of periodic points. Then the
 intersection of $W^s(x)$ and $W^u(y)$ at $z$ is transversal.
 The same argument of last paragraph then shows that $z\in\La$.
 
 \end{proof}

 \begin{Lemma}\label{closedint}\quad
 
 Let $x$ be an interior point of the homoclinic class $\La$.
 Then there is $\e=\e(x)>0$ such that the set  $Q\cap \ov{B_\e(x)}$
 is closed,  where $B_\e(x)$ is
 the $\e$-ball centered at $x$.
 
 \end{Lemma}
 
 \begin{proof}
 Observe that for $\de>0$ small enough, the set
 $$
 Q_\de :=\bigcup_{x\in\La}\ov{W^{ss}_\de(x)}\cup\ov{W^{uu}_\de(x)}
 $$
 is closed. It is enough to prove that for $\e=\e(x)>0$ sufficiently small,
 $B_\e(x)\cap Q=B_\e(x)\cap Q_\de$. 
 
 Let $\cR$ be a small surface transversal
 to the vector field of $\phi_t$ containing $x$ in its interior, which
 is a topological rectangle whose sides are in two weak stable
 manifolds and two weak unstable manifolds:
   $$
   \partial \cR \subset  \bigcup_{i=1,2} W_\e^s(x_i)\cup W_\e^u(x_i),
    \qquad\text{ with } x_1,\;x_2\in\La.
   $$ Since $x$ is an interior
 point of $\La$, the diameter of $\cR$ can be chosen arbitrarily small.
 Let $\tau>0$ be such that $\phi_t(\cR)\cap\phi_s(\cR)=\emptyset$ for
 all $-\tau\le t < s\le \tau$. For simplicity assume that $\tau=2$.

 Recall that (cf. \cite[\S 6]{HPPS} or \cite{HPS}) 
 since $\La$ is uniformly hyperbolic, its local
 weak invariant manifolds $W^s(x)$, $W^u(x)$, $x\in\La$ are
 graphs of functions 
 $A(x):E^s(x)\oplus\langle \tfrac d{dt}\phi_t\rangle\to E^u(x)$
 whose domain contains a ball of some fixed radius $2\,\rho_1>0$.
 In particular their diameters are larger than $\rho_1$.
 Similarly for the strong invariant manifolds.
 
 We take the diameter of $\cR$ small enough so that if
 $z\in\phi_{s}(\cR)\cap\La$, $s\in[-1,1]$, then 
 $\big[\cup_{t\in[-\de,\de]}W^{ss}_\de(\phi_t(z))\big]\cap\phi_{s}(\cR)$ 
 and  $\big[\cup_{t\in[-\de,\de]}W^{uu}_\de(\phi_t(z))\big]\cap\phi_{s}(\cR)$
 are differentiable curves which cross the rectangle
 $\phi_{s}(\cR)$. Now take $\e=\e(x)>0$ such that
 $\ov{B_\e(x)}\subset \phi_{[-1,1]}(\cR)$.

 \begin{figure}
 \includegraphics[scale=.5]{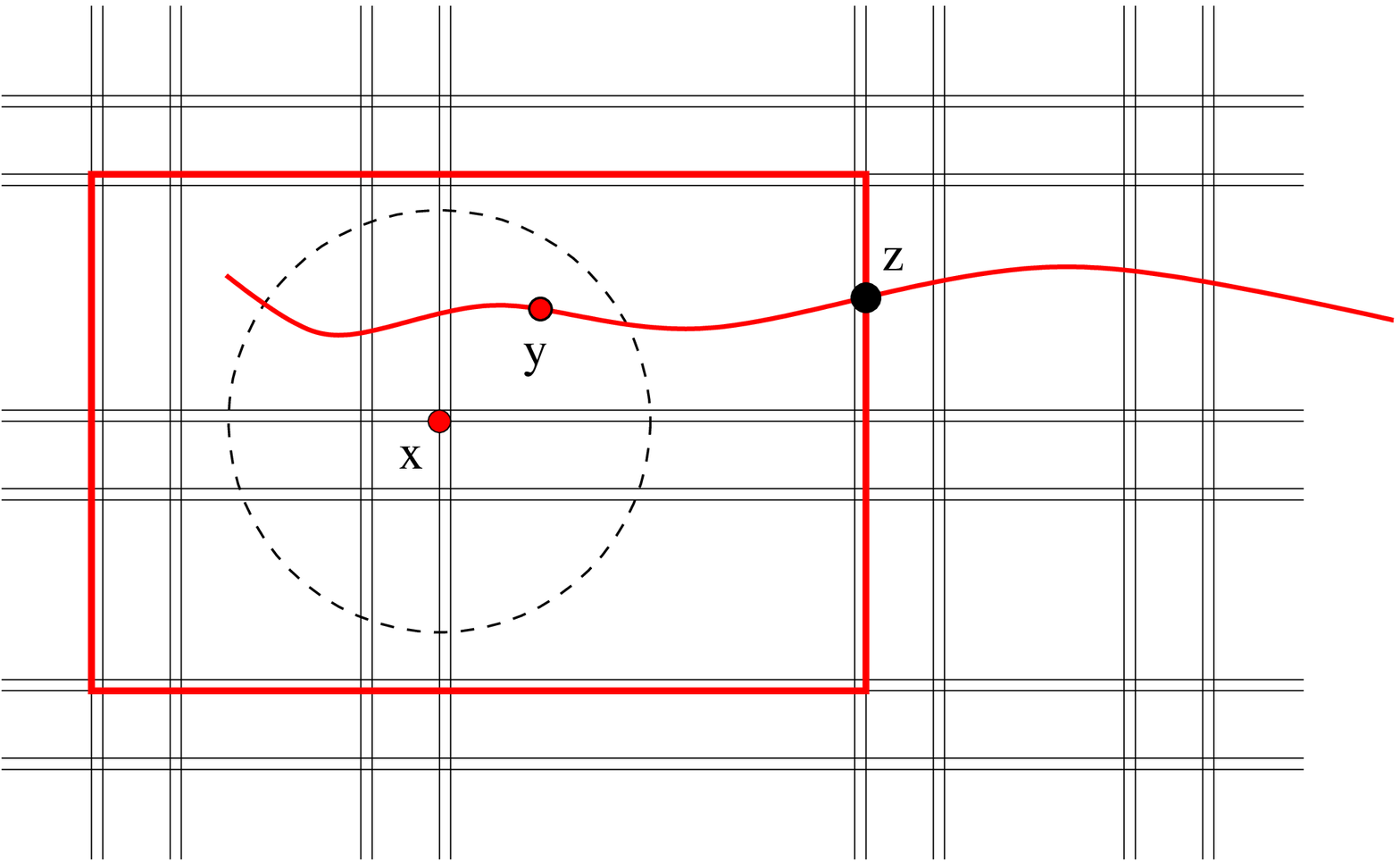}
 %\refstepcounter{figure}\label{rect}
 \caption{}\label{rect}
 \end{figure}

 Let $y\in Q\cap \phi_{s}(\cR)$, $s\in[-1,1]$, and suppose that $y\in
 W^{s}(w)$, with $w\in\La$. The case $y\in W^{u}(w)$, $w\in\La$ is 
 similar. Consider the connected component $\Ga$ of the intersection
 of the weak stable manifold  $W^s(w)$ with the rectangle $\phi_{s}(\cR)$. 
  Suppose first that $\Ga$ intersects the boundary of $\phi_{s}(\cR)$,
  necessarily at one of its endpoints $z$, as shown in
  figure~\ref{rect}. 
  Then $z$ is in the intersection of $W^s(w)$
  with $W^u(x_1)$ or $W^u(x_2)$.
   Therefore, by lemma~\ref{homclass},
    $z$ is in the homoclinic class $\La$.
   By our choice of the diameter of $\cR$, we have that
   $y\in \cup_{t\in[-\de,\de]}W^{ss}_\de(\phi_t(z))$.
   Hence $y\in Q_\de$.
   
   Now suppose that $\Ga$ does not intersect the boundary of  
   $\phi_{s}(\cR)$. The curve $\Ga=W^s(w)\cap
   \phi_s(\cR)\subset\inte(\cR)$ must have at least one well
   defined endpoint $z$, which is in the orbit of some $w\in\La$.
   By our choice of the size of $\cR$, the stable manifold
   $W^s(z)=W^s(w)$ of $z\in\La\cap\phi_s(\cR)$ crosses $\phi_s(\cR)$.
   Therefore $\Ga\cap\partial \,[\phi_s(\cR)]\ne\emptyset$ and 
   hence this case  can not occur.
   
   We have proved that $Q\cap \ov{B_\e(x)}\subseteq
   Q\cap\phi_{[-1,1]}(\cR)\subseteq Q_\de$.

\end{proof}

 % Let $\gamma$ be a simple closed geodesic. For example, 
 % it can be obtained by a minimax method (see Birkhoff~\cite{Birk2} or  
 % e.g. Struwe~\cite{struwe0}). Then $\ga$ is the boundary of an embedded disc
 % $D$ on $S^2$. Let $\Sigma$ be the set of unit tangent vectors on 
 % $S^2$ based at the points of $\ga$ and pointing inwards $D$.
 % Since $\ga$ is a geodesic, the surface $\Sigma$ is transverse 
 % to the geodesic flow. Birkhoff~\cite[page...]{Birk2} shows that if $g_0$ 
 % has positive curvature, then $\Sigma$ intersects all the orbits 
 % of the geodesic flow  infinitely many times, 
 % except for the orbit  $(\ga,\dot{\ga})$ and its flip
 % $(\ga,-\dot\ga)$.
 % It follows that the Poincar\'e's return times to $\Sigma$ are finite 
 % (but possibly not bounded) and that the Liouville's measure 
 % induces a smooth probability measure $\mu$, invariant 
 % under the return map $f:\Sigma\to\Sigma$.

 Since $g$ is a Kupka-Smale metric on $S^2$, by 
 proposition~\ref{conformal} we can lift the geodesic flow 
 to a reparametrization of the Reeb flow of a non-degenerate
 tight contact form on $S^3$.

 \medskip
 
 \noindent
 {\bf The dynamically convex case.}
 
 \smallskip
 
 In this case there is a closed orbit $\ga$
 for the lifted geodesic flow $\psi_t$ to $S^3$ and
 a smooth 2-dimensional open disk $\di\subset S^3$ such that
 $\partial\di=\ga$, $\di$ is transversal to the 
 lifted geodesic flow $\psi_t$ and  intersects
 every orbit in $S^3\setminus\ga$ infinitely many times.
 Moreover, the differential $d\la$ of the lifted 
 Liouville's form $\la$ is non-degenerate on
 $\di$, the total area $\int_\di d\la=\text{period}(\ga)$ 
 is finite and the Poincar\'e's return map
 $f:\di\to\di$ is $C^2$ and preserves $\la$.
 Also the first return time $\tau:\di\to ]0,+\infty[$
 is uniformly bounded above and below.
 
 Let $\hLa$ be the lift of the basic set $\La$ to $S^3$
 and let $\Ga:=\hLa\cap\di$.
 Let $\mu$ be the normalized $d\la$-measure on $\di$.
 The invariant measure for $\psi_t$ induced by
 $\mu$ is the lift of the Liouville's measure $m$.
 Then $\mu(\Ga)=0$.
 Since there are $0<a<b$ such that the return time $\tau$ 
 is uniformly bounded: $a<\tau(x)<b$, $\forall x\in\di$,
  the set $\Ga$ is a uniformly hyperbolic set for 
 the Poincar\'e map $f(x)=\psi_{\tau(x)}(x)$.
 
 The hyperbolic set $\Ga$ may contain the boundary
 closed orbit $\ga=\partial\di$. In that case $\Ga$ 
 accumulates on $\ga$ and in particular is not closed.
 But nevertheless it is locally compact {\it inside} $\di$.
 
  Since $\Ga$ is non-trivial it has an interior point $x$.
   Let $\epsilon:=\min\{\e(x),\,\tfrac 12\,d(x,\partial
  \di)\}>0$, where $\e(x)$ is from lemma~\ref{closedint}.
  Let $B_\epsilon(x)$ be the ball of radius $\epsilon$ in $\di$.
   Let $\cQ$ be the lift of $Q$ intersected with $\D$.
  Then $\cQ=W^s(\Ga)\cap W^u(\Ga)$.
  Since the points in $\cQ\setminus \Ga$ are wandering,
  $\mu(\cQ)=\mu(\Ga)+\mu(\cQ\setminus\Ga)=0$.
  Let $R$ be a connected component of $\di\setminus Q$
  which is contained in $B_\e(x)$.
  Such $R$ always exist because $x$ is an interiror point
  of $\Ga$, $\ov{B_\e(x)}\cap\cQ$ is closed by lemma~\ref{closedint}
  and $\mu(\cQ)=0$.
   % since $\Ga$ is a non-trivial basic set, it has no isolated points.
 Moreover $R$ is a topological rectangle whose sides are
 pieces of stable and unstable manifolds of points in $\Ga$:
   $$
  \partial R \subset W_\e^s(x_1)\cup W_\e^s(x_2)
  \cup W_\e^u(y_1)\cup W_\e^u(y_2). 
  $$
 Since $R$ has positive $d\la$-measure, by Poincar\'e's recurrence
 there exists $N>0$ such that $f^N(R)\cap R\ne \nulo$. Since $Q$ is 
 invariant and $R$ is a connected component of $\di\setminus Q$,
 we have that $f^N(R)=R$. Since\footnote{Another proof from this point
 is the following: since $R$ is homeomorphic to  a disc,
 and $f^N$ preserves a finite measure,
 by Brouwer's translation theorem there is a fixed point $y$ for $f^N$
 inside $R$. The point $y$ belongs to a periodic orbit for $\psi$.
 Since $\ov{\Per(\psi)}$ is uniformly hyperbolic, if $\epsilon>0$ is small
 enough then $\emptyset\ne W^s(y)\cap\partial R\subset W^s(y)\cap W^u(\La)$.
 Then $y$ is in the homoclinic class $\La$. Since $y\in R$,
 this contradicts the choice of $R$. We shall use this argument in the 
 non-dynamically convex case. We present another argument here to show that
 Brouwer's translation theorem is not needed for the dynamically 
 convex case or for the positively curvature case.}
  $f^N$ is continuous,  $f^N(\partial R)=\partial R$.
 Since local stable manifolds are sent to local stable manifolds,
 using twice the period if necessary, we have that 
 $f^{2N}(\ov{I_1})\subset \ov{I_1}$, where
 $I_1:=W^s_\e(x_1)\cap R$.
 By the intermediate value theorem there is a fixed point 
 $y\in\ov{I_1}$, $f^{2N}(y)=y$. Since $g\in\cF^1(S^2)$, 
 $y$ is a hyperbolic point. Since $I_1$ is an invariant curve 
 under $f^{2N}$ and $y\in\ov{I_1}$, we have that $I_1\subset W^s(y)$.
 Parallel sides of the rectangle $R$ are disjoint
 $W_\e^u(x_1)\cap W^u_\e(x_2)=\nulo$. We can assume that
 $y\notin W_\e^u(x_1)$.
 Since $g$ is Kupka-Smale, $W_\e^u(x_1)$ intersects transversally
 $I_1$ at one of its endpoints $w$. Since $w\in I_1\subset W^s(y)$,
 we have that $f^{2Nk}(W^u_\e(x_1))$ is transversal to $I_1$ and
 contains the point $f^{2Nk}(w)$ which approaches $y$ as $k\to+\infty$.
 Then  $\nulo\ne f^{2Nk}(W^u_\e(x_1))\cap\inte(R) \subset Q\cap\inte(R)$ 
 for some $k\in\na$. But this contradicts the choice of $R$.
 
 \qed

 \bigskip
 
 \noindent{\bf The non-dynamically convex case.}
 
 \bigskip

 Let $\psi_t$ be the lift of the geodesic flow to $S^3$.
 Then $\psi_t$ is the Reeb flow of a tight contact form
 and it is Kupka-Smale. The lift $\oPs$ of $\oPg$ is also a 
 hyperbolic set. Let $\hLa$ be the lift of $\La$ and $\hQ$ the lift of
 $Q$.
  
 For $\ep>0$ and $x\in\oPs$, let $W^s_\epsilon(x)$, $W^u_\ep(x)$
 be the $\ep$-balls in $W^s(x)$, $W^u(x)$.
 Since the set $\oPs$ is hyperbolic (cf.~\cite[p. 6.4.13]{HK},
 \cite{Bowen6}), 
 for all $\ep>0$ there
 exists $\de=\de[\ep]>0$ such that 
 \begin{equation}\label{deba}
 \text{ if }\quad y,\,z\in\oPs\quad\text{ and }\quad
 d(y,z)<\de[\epsilon]\quad\text{ then }\quad 
 W_\ep^u(y)\cap W_\ep^s(z)\ne \nulo.
 \end{equation}

 Let $\Si$ be a rigid surface such that $\Si\cap \hLa\ne \nulo$.
 Let $x\in\hLa\cap\Si$ be an interior point of $\hLa\cap\Si$.  Let
 $$
 \e=\min\left\{\de\left[\tfrac 12 \,d(x,\partial\Si)\right],\,
 A,\,\e(x)\,\right\}>0,
 $$ 
 where $A>0$ is from proposition~\ref{diamF}, $\e(x)>0$ is from 
 lemma~\ref{closedint} and $\de[\,\cdot\,]$ is a function 
 satisfying~\eqref{deba}.
 
  Let $B_\e(x)$ be the ball of radius $\e$ in $\Si$
 centered at $x$. By lemma~\ref{closedint},
 $\ov{B_\e(x)}\cap\hQ\cap\Si$ is closed. 
 The contact 1-form $\la$ is the lift of the Liouville's form
 on $T^1S^2$. The form $d\la$ on $\Si$ is finite and is
 preserved by the return  map to $\Si$ it induces a smooth
 invariant probability $\mu$ on $\Si$ . 
 Since the flow $\psi_t$ is not Anosov,
 $\mu(\hLa\cap\Si)=0$. Since the points in 
 $(\hQ\cap\Si)\setminus \hLa$ are wandering, by Poincar\'e's
 recurrence, $\mu[(\hQ\cap\Si)\setminus\hLa]=0$.
 Thus $\mu(\hQ\cap\Si)=0$.
 Therefore $\ov{B_\e(x)}\setminus\hQ$ is an open non-empty set.

 Let $R$ be a connected component of $B_\e(x)\setminus \hQ$
 whose closure is in $\inte(B_\e(x))$. Such component exists
 because $x$ is not isolated in $\hLa$. Actually, $R$ is an open
  rectangle with two sides in each of $W^s(\hLa)$ and $W^u(\hLa)$.
 Let $\tau_n:\Si\to [0,+\infty[\cup\{+\infty\}$ and 
 $F_n:\Si\to\Si$ be the return times and return maps to $\Si$,
 i.e.
 \begin{align*}
 \tau_0\equiv 0,
 \qquad
 F_n(x)=\psi_{\tau_n(x)}(x), 
 \qquad
 \tau_{n+1}(x)=\inf\,\{\,t>\tau_n(x)\;|\; \psi_t(x)\in\Si\,\}.
 \end{align*}
  By Poincar\'e's recurrence, there is some $N>0$ such that
  $F_N(R\cap[\tau_N<+\infty])\cap R\ne \nulo$.
  
  Suppose that\footnote{Observe that a boundary segment of $R$ could still
  be in a stable manifold of a binding orbit and never return to
  $\Si$.} $\inte(R)=R\subset [\tau_N<+\infty]$.
  Then $F_N\vert_R$ is a local diffeomorphism and then $F_N(R)$ is connected.
  If $F_N(R)\cap \partial R\ne \nulo$ then there is a point
  $q\in\partial R\cap F_N(R)\subset \hQ\cap F_N(R)$. Since the set $\hQ$
  is  invariant,  there is $q'\in R\cap \hQ$, which contradicts
  the choice of $R$. Hence $F_N(R)\subseteq R$. Since the map
  $F_N:R\to R$ is area preserving and the $d\la$-area of $R$ is finite, 
  it can not have a translation  domain. 
  By Brouwer's translation theorem~\cite{Fra3}, $F_N$ has a
  fixed point in $R$. Then the flow $\psi_t$ has a periodic point
  $p\in R\cap \oPs$.  By the definition of $\de=\de[\epsilon]$ above, we have that
  $W^s(p)\cap W^u(x)\ne 0$ and $W^u(p)\cap W^s(x)\ne 0$.
  Then $p$ is in the same homoclinic class as $x$, thus $p\in\hLa$.
  But this contradicts the choice of $R$. 
  
  Then $R\cap[\tau_N=+\infty]\ne\nulo$.
  By the choice of $N$ there exists $q\in R\cap[\tau_N<+\infty]$ such
  that $F_N(q)\in R$. Let $B$ be the connected component of 
  $R\cap[\tau_N<+\infty]$ containing $q$. Since $[\tau_N<+\infty]$ is
  open, we have that $B$ is open. Since $R\cap[\tau_N=+\infty]\ne\nulo$,
   $\tau_N$ is unbounded in $B$. By proposition~\ref{diamF},
  $\diam F_N(B)>A$. Since $\diam(R)<\e\le A$, we have that
  $F_N(B)\not\subset R$. But since $F_N(B)$ is
  connected, the same argument as above shows that $F_N(B)\subseteq
  R$ because $F_N(B)\cap\partial R\subseteq F_N(B)\cap\hQ=\emptyset$. This is a contradiction.
  
  \qed

  \def\cprime{$'$} \def\cprime{$'$} \def\cprime{$'$}
 \providecommand{\bysame}{\leavevmode\hbox to3em{\hrulefill}\thinspace}
 \providecommand{\MR}{\relax\ifhmode\unskip\space\fi MR }
 \providecommand{\MRhref}[2]{%
   \href{http://www.ams.org/mathscinet-getitem?mr=#1}{#2}
 }
 \providecommand{\href}[2]{#2}

\end{document}